\documentclass[onefignum,onetabnum]{siamonline190516}

\usepackage{lipsum}
\usepackage{amsfonts,amsmath}
\usepackage{graphicx}
\usepackage{epstopdf}
\usepackage{algorithm}
\usepackage{algpseudocode}
\algdef{SE}[SUBALG]{Indent}{EndIndent}{}{\algorithmicend\ }%
\algtext*{Indent}
\algtext*{EndIndent}
\usepackage{wrapfig}
\usepackage{tikz-cd}
\ifpdf
  \DeclareGraphicsExtensions{.eps,.pdf,.png,.jpg}
\else
  \DeclareGraphicsExtensions{.eps}
\fi

\usepackage{
    mathtools,      
    amssymb,        
    amsfonts,       %
    thmtools,       %
    mathrsfs,       
}

\def\N{\mathbb N}

\def\R{\mathbb R}

\def\E{\mathbb E}

\def\mC{\mathcal{C}}
\def\tC{\tilde{\mathcal{C}}}

\def\b{\beta}

\def\l{\lambda}

\def\ben*{\begin{eqnarray*}}
\def\eeqn*{\end{eqnarray*}}

\usepackage{cancel}

\definecolor{PenGreen}{RGB}{37,158,108}
\definecolor{cornflower}{RGB}{95,148,201}
\definecolor{blue(ncs)}{rgb}{0.0, 0.53, 0.74}

\usepackage{
    nameref,         
}

\usepackage{subfig}

\theoremstyle{definition}

\usepackage{enumitem}
\setlist[enumerate]{leftmargin=.5in}
\setlist[itemize]{leftmargin=.5in}

\newsiamremark{remark}{Remark}
\newsiamremark{hypothesis}{Hypothesis}
\crefname{hypothesis}{Hypothesis}{Hypotheses}
\newsiamthm{claim}{Claim}

\headers{Higher Order Unscented Transform}{D. Easley and T. Berry}

\title{A Higher Order Unscented Transform \thanks{Submitted to the editors DATE.
\funding{Dr. Berry was supported by NSF under grant no.~1854204.}}}

\author{Deanna Easley \thanks{Dept. of Mathematical Sciences, George Mason University, Fairfax, VA 
  (\email{deasley2@gmu.edu}).}
\and Tyrus Berry\thanks{Dept. of Mathematical Sciences, George Mason University, Fairfax, VA 
  (\email{tberry@gmu.edu}, \url{http://math.gmu.edu/\string~berry/}).}
}

\usepackage{amsopn}

\ifpdf
\hypersetup{
  pdftitle={A Higher Order Unscented Transform},
  pdfauthor={D. Easley and T. Berry}
}
\fi

\begin{document}

\maketitle

\begin{abstract}
We develop a new approach for estimating the expected values of nonlinear functions applied to multivariate random variables with arbitrary distributions. Rather than assuming an particular distribution, we assume that we are only given the first four moments of the distribution. The goal is to summarize the distribution using a small number of quadrature nodes which are called $\sigma$-points.  We achieve this by choosing nodes and weights in order to match the specified moments of the distribution. The classical scaled unscented transform (SUT) matches the mean and covariance of a distribution. In this paper, introduce the higher order unscented transform (HOUT) which also matches any given skewness and kurtosis tensors. It turns out that the key to matching the higher moments is the rank-1 tensor decomposition. While the minimal rank-1 decomposition is NP-complete, we present a practical algorithm for computing a non-minimal rank-1 decomposition and prove convergence in linear time. We then show how to combine the rank-1 decompositions of the moments in order to form the $\sigma$-points and weights of the HOUT.  By passing the $\sigma$-points through a nonlinear function and applying our quadrature rule we can estimate the moments of the output distribution. We prove that the HOUT is exact on arbitrary polynomials up to fourth order. Finally, we numerically compare the HOUT to the SUT on nonlinear functions applied to non-Gaussian random variables including an application to forecasting and uncertainty quantification for chaotic dynamics.
\end{abstract}



\section{Introduction}
A fundamental problem in uncertainty quantification is to approximate the expectation a function $f$ applied to a random variable $X$ sampled from a probability measure $dp$, namely
\begin{equation}\label{integral} \mathbb{E}[f(X)] = \int f(x)\, dp. \end{equation} 
Even when everything is known this can be a challenging computation in high dimensions, and the problem is often compounded by uncertain or incomplete knowledge of $f$ and $dp$.  Moreover, in most problems of interest $f$ has an extremely complex form, for example $f$ may encapsulate the solution of a differential equation and the computation of some feature of interest on the solution.  So we may not be able to assume that $f$ is known in an explicit form, but instead that $f$ or an approximation to $f$ is available only as a black-box computational scheme which can take inputs $x$ and produce outputs $f(x)$.  Similarly, the type of partial knowledge of the probability measure can vary widely, we may have an explicit expression for a density function $p(x) = dp/dx$ (if it even exists), or we may only have some samples of $dp$ or estimates of some of the moments. 

In any of the above circumstances, the problem of approximating \eqref{integral} can be approached as a problem of numerical quadrature (also known as cubature when $x$ has dimensionality greater than one, we will use the term quadrature for both).  A quadrature is an approximation of the form
\begin{equation} \label{approx} \mathbb{E}[f(X)]\approx \sum_{i=1}^N w_i f(x_i) \end{equation}
where $x_i$ are called nodes and $w_i$ are called weights.  The goal is to find a small number of nodes and weights that accurately represent the probability measure for a large space of functions $f \in \mC$.  A common strategy in quadrature methods is to choose nodes and weights so that the above approximation is actually an equality for all $f$ in some finite dimensional subspace $\tC \subset\mC$ (such as a space of polynomials up to a fixed degree).  For $f$ outside of $\tC$ we can then attempt to bound the error in the approximation \eqref{approx} if we can control the error between $f$ and its projection into $\tC$.   When $f$ is sufficiently smooth and $dp$ is concentrated in a small region, then it is reasonable to approximate $f$ using the space of polynomials up to a fixed degree.  In this case, ensuring that \eqref{approx} holds with equality for all polynomials up to degree $k$ leads to the so-called moment equations,
\begin{equation} \label{moments} m_{j_1,...,j_n} = \mathbb{E}\left[X_{1}^{j_1}X_{2}^{j_2}\cdots X_{n}^{j_n}\right] = \sum_{i=1}^N w_i x_1^{j_1}x_2^{j_2} \cdots x_n^{j_n} \end{equation}
for all $j_1 + j_2 + \cdots j_n \leq k$.  In other words we are asking that the empirical moments of the nodes $x_i$, weighted by discrete probabilities $w_i$, exactly agree with the true moments $m_{j_1,...,j_n}$ of the distribution.  When $k=2$ the moment equations specify that weighted nodes must match the mean vector and covariance matrix of the true distribution, and this is achieved with the so-called \emph{Scaled Unscented Ensemble} (SUT) \cite{julier} (see Section \ref{SUT} for details).  

The quadrature approach is an alternative to stochastic quadrature methods such as Monte-Carlo quadrature as is commonly used in particle filtering. Stochastic quadratures use random variables $X_i$ to build quadrature rules such that
 \begin{equation} \label{approxMC} \mathbb{E}[f(X)] \approx \mathbb{E}\left[\sum_{i=1}^N w_i f(X_i)\right] \end{equation}
however, the computed value $\sum_{i=1}^N w_i f(X_i)$ will be stochastic.  This means that in addition to possible approximation error in \eqref{approxMC}, we also have an error due to the variance of the random variable  $\sum_{i=1}^N w_i f(X_i)$.  While it is often easier to design stochastic quadrature methods where the approximation error in \eqref{approxMC} is small or even zero, for many problems controlling the variance error requires a large number of random variables $X_i$ and hence a large number of function evaluations.  When $f$ is very complex is can be beneficial to have a small deterministic ensemble and accept the quadrature error in \eqref{approx} in order to avoid the variance error of a stochastic quadrature.

The problem \eqref{integral} is often part of a larger problem such as filtering \cite{filter1}, particle filtering \cite{filter2}, adaptive filtering \cite{berry2013adaptive}, smoothing \cite{smoother1}, parameter estimation \cite{param1,param2,param3} and even model-free filtering \cite{modelfree}.  In all these applications it can be beneficial to have deterministic approximation of \eqref{integral} to improve the stability of the overall algorithm. For example, filters built on random ensembles can fail catastrophically since they generate many such ensembles can generate realizations that would normally have very low probability but lead to perverse behavior \cite{harlim2010catastrophic,anderson2007adaptive}.  Similarly, gradient based optimization method for parameter estimation will need to carefully account for any stochasticity in the objective function, so replacing a stochastic quadrature with a deterministic quadrature can be desirable in certain applications. 

The highly successful Uscented Kalman Filter (UKF) \cite{filter1} is based on the SUT along with the many other algorithms mentioned above. A closely related technique called Cubature Kalman Filters (CKF) \cite{cubaturefilter} follow a similar strategy and are typically designed to achieve a high degree of exactness under a Gaussian assumption on the distribution.  In this paper we avoid any assumptions on the distribution and instead achieve a high degree of exactness based purely on moment matching. Whereas the UKF (and implicitly most CKFs) only require the rank-1 decomposition of the covariance matrix, the Higher Order Unscented Transform (HOUT) requires the rank-1 decomposition of higher order tensors such as the skewness and kurtosis.  
The rank-1 tensor decomposition of a $k$-tensor is defined by vectors $v_i$ such that, \begin{equation}\label{r1decomp} T=\sum_{i=1}^p v_i^{\otimes k}.\end{equation} 
The minimal value of $p$ such that the above decomposition exists is the called the rank of $T$.  

Ideally, we would like an exact rank-1 tensor decomposition \eqref{r1decomp} with the minimum possible number of vectors, however this turns out to be an NP-complete problem \cite{NP1,NP2}.  Instead, we will use an effective algorithm for obtaining an approximate rank-1 decomposition up to an arbitrary tolerance.  The algorithm was originally suggested by \cite{EigApprox}, and it works by repeatedly subtracting the best rank-1 approximation to a tensor until the norm of the residual is less than any desired tolerance.  Many methods have been developed based on this idea (see \cite{survey} and citations therein) and in \cite{convergence} it was proven to converge but without any convergence rate.  In Section \ref{approxrank1} we give the first proof that this algorithm converges linearly and we derive an upper bound on the convergence rate.  While the approximate decomposition typically requires many more vectors than the minimal rank-1 decomposition, it avoids the NP-completeness of that problem and gives us an effective algorithm.  

In Section \ref{sec:background} we briefly review the SUT and some tensor facts and notation including the Higher Order Power Method (HOPM) \cite{HOPMoriginal} that we will use for finding tensor eigenvectors.  Based on the HOPM, we prove the convergence of the approximate rank-1 decomposition algorithm in Section \ref{approxrank1}.  This proof also requires some new inequalities relating the maximum eigenvalue of a tensor to the entries of the tensor, and these inequalities are likely to be of independent interest.  In Section \ref{sec:HOUT} we introduce the Higher Order Unscented Transform (HOUT) which generalizes the SUT in order to match arbitrary skewness and kurtosis tensors. The HOUT gives a quadrature rule with degree of exactness four that is applicable to arbitrary distributions.  Finally, we demonstrate the HOUT on various non-Gaussian multivariate random variables on complex nonlinear transformations in Section \ref{sec:numerics} and briefly conclude in Section \ref{sec:conclusions}.

\section{Background}
\label{sec:background}

We start by reviewing the Scaled Unscented Transform (SUT) in Section \ref{SUT} which has degree of exactness two.  We then briefly introduce our tensor notation in Section \ref{tensors} and tensor-vector products and tensor norms in Section \ref{tensornorms}.  Finally, in Section \ref{tensoreigs} we review tensor eigenvectors and eigenvalues and the Higher Order Power Method (HOPM) \cite{HOPMoriginal} for finding them.

\subsection{Scaled Unscented Transform}\label{SUT}

The Scaled Unscented Transform (SUT) was introduced by Julier and Uhlmann in \cite{julier} and further developed in \cite{julier1995,julier1997,julier1998skewed,julier2002scaled,julier2004unscented}. The fundamental goal of this paper is to generalize their method to higher order moments.  This work was started in \cite{julier1998skewed} which worked on matching the skewness, and below we show that rank-1 tensor decompositions are the key to generalizing their approach.  

The SUT uses the mean and covariance of a distribution to choose quadrature nodes and weights such that the quadrature rule has degree of exactness 2.  Degree of exactness $k$ means that a quadrature rule is exact for computing the expectation of polynomials up to degree $k$.  The fundamental insight of Julier and Uhlmann is that achieving degree of exactness 2 is equivalent to matching the first two moments of the distribution. Moreover, they showed that this can be efficiently accomplished using a matrix square root of the covariance matrix (which is a rank-1 decomposition of the covariance matrix). 
\vspace{0.75em}

\begin{definition}[$i$th column of the symmetric matrix square root of $A$]
Let $A$ be a $d\times d$ matrix. We define the {\boldmath {\bf $i$th column of the symmetric matrix square root of $A$}}, denoted $\sqrt{A}_i$, by \[\sum_{i=1}^{d} \sqrt{A}_i^{\otimes 2}=\sum_{i=1}^d \sqrt{A}_i\sqrt{A}_i^\top=A.\]
\end{definition}
The notation $v^{\otimes k}$ will be defined below. Note that the following definition can use any matrix square root but we have found empirically that the unique symmetric matrix square root has the best performance.
\begin{definition}[The Scaled Unscented Transform (SUT) \cite{julier}]
Let $dp$ be a probability measure with mean $\mu\in\R^d$ and the covariance $C\in\R^{d\times d}$. Then for some $\beta\in\R$ the 
{\boldmath{\bf$\sigma$--points}} are defined by 
\[
\sigma_i = \begin{cases} 
 \mu  & \hbox{if } \  i=0 \\
     \mu+\beta\sqrt{C}_i & \hbox{if } \  i=1,\dots, d \\
    \mu-\beta\sqrt{C}_{i-d} & \hbox{if } \  i=d+1,\dots, 2d 
   \end{cases}\]
and the corresponding {\boldmath {\bf weights}} are defined by
\[w_i = 
\begin{cases}
    1-\frac{d}{\beta^2} & \hbox{if } \  i=0 \\
    \frac{1}{2\beta^2} & \hbox{if } \  i=1,\dots, 2d \\
   \end{cases}
\]
\end{definition}
We note that the choice of $\beta$ can have significant impact on the effectiveness of the transform.

\begin{remark}
The absolute condition number of the Scaled Unscented Transform is \\bounded above by \[\sum_{i=0}^{2d}|w_i|=\left|1-\frac{d}{\beta^2}\right|+\frac{d}{\beta^2}.\] If $\beta\ge\sqrt{d}$, then \(\displaystyle\sum_{i=0}^{2d}|w_i|=1.\) If $\beta<\sqrt{d}$, then \(\displaystyle\sum_{i=0}^{2d}|w_i|=\frac{2d}{\beta^2}-1.\)
\end{remark}

\begin{theorem}[Empirical mean and Empirical covariance \cite{julier}]
For an arbitrary $\beta$, we have
\[\mu = \E[X] = \sum_{i=0}^{2d} w_i\sigma_i \hspace{20pt}{\rm and}\hspace{20pt} C = \E[(X-\mu)(X-\mu)^\top] = \sum_{i=0}^{2d} w_i(\sigma_i-\mu)(\sigma_i-\mu)^\top 
\]
and if $q:\mathbb{R}^d \to \mathbb{R}$ is a polynomial of degree at most $2$, we have, $\mathbb{E}[q(X)] = \sum_{i=0}^{2d} w_i q(\sigma_i)$.
\end{theorem}
We should note that if the distribution has zero skewness, such as a Gaussian distribution, then the symmetry of the nodes yields degree of exactness $3$, and, in the specific case of a Gaussian distribution the choice $\beta = \sqrt{3}$ acheives degree of exactness $4$ \cite{julier,julier2002scaled,julier2004unscented}. The choice $\beta = \sqrt{d}$ is often called the \emph{unscented transform} and sets $w_0=0$ so that only $2d$ of the $\sigma$-points are required.  The ability of the SUT to match the first four moments of the Gaussian distribution has led some to associate the SUT with a Gaussian assumption, however this is not required and degree of exactness $2$ is acheived for arbitrary distributions.  Our goal is to generalize the unscented transform to higher moments, which are tensors.

\subsection{Tensors}\label{tensors}
Tensors are essentially multidimensional matrices, which will be used to conveniently express the notions of covariance, skewness and kurtosis in a similar fashion.

\begin{definition}[$k$-order tensor] For positive integers $d$ and $k$, a
 tensor $T$ belonging to $\R^{d^k}$ is called a {\boldmath\bf $k$-order tensor} or simply a {\boldmath\bf$k$-tensor}. 
\end{definition}

In particular, a vector in $\R^d$ can be viewed as a first order tensor and a $d\times d$ matrix as a second order tensor. Let $x\in\R^d$.  We note that the outer product $xx^\top$ yields a $d\times d$ matrix whose $ij$-entry can be represented as 
\[(xx^\top)_{ij}=x_ix_j = (x\otimes x)_{ij} = (x^{\otimes2})_{ij}.\]
We generalize this process of forming higher order tensors from vectors with following definition.
\begin{definition}[$k$th-order tensor product]\label{tensorprod}
Let $v\in\R^d$ and $k$ be a positive integer then the {\boldmath\bf$k$th-order tensor product} is a $k$-tensor denoted \[v^{\otimes k}=\underbrace{v\otimes v\otimes\cdots\otimes v}_{k \ \rm{times}},\] the elements are given by $(v^{\otimes k})_{i_1,\dots,i_k} = v_{i_1}\dots v_{i_k}$.
\end{definition}
Definition \ref{tensorprod} immediately connects tensor products to the moments of a distribution since we can represent the covariance as $C = \E[(X-\mu)^{\otimes2}] = \int (x-\mu)^{\otimes2}\, dp(x)$ so that the skewness $S$ and kurtosis $K$ can be defined as 
\[ S = \int (x-\mu)^{\otimes3} \, dp(x) \hspace{30pt} K = \int (x-\mu)^{\otimes4} \, dp(x), \]
so that, for example, 
\[ S_{ijk} =\int \left((x-\mu)^{\otimes3}\right)_{ijk} \, dp(x) = \int (x-\mu)_i(x-\mu)_j(x-\mu)_k \, dp(x). \]
The following definition generalizes the notion of a rank-1 matrix to tensors. We will return to the general notion of tensor rank in a later section.
\begin{definition}[Rank-1 Tensor]
Let $T\in\R^{d^k}$ then $T$ is called a {\bf rank-1 tensor} if there exists a $v\in\R^d$ such that \[v^{\otimes k}=T.\]
\end{definition}
For tensors that are not rank-1, one may seek a rank-1 decomposition of the form,
\[ T = \sum_{\ell=1}^p v_\ell^{\otimes k} \]
and the minimum $p$ for which such a decomposition exists is called the \emph{rank} of the tensor $T$.  This notion of rank agrees with the classical notion of matrix rank in the case of second order tensors but many of the properties of matrix rank do not generalize to higher order tensors \cite{tensor1,NP1,NP2,subtract1,subtract2}.

\subsection{Tensor Multiplication and Tensor Norms}\label{tensornorms}

To discuss how tensor multiplication works, let us first look at the simplest case where we multiply a 2-tensor with a 1-tensor. Recall that for a matrix $A\in\R^{d\times d}$ and $v\in\R^d$, the matrix-vector multiplication $Av$ is given by $(Av)_i = \sum_{j=1}^d A_{ij}v_j$ so we define two natural tensor-vector products
\ben*
(A\times_1v)_i &=& \sum_{j=1}^d A_{ji}v_j = (A^\top v)_i \hspace{20pt}{\rm and}\hspace{20pt}
(A\times_2v)_i = \sum_{j=1}^d A_{ij}v_j = (A v)_i
\eeqn*
Analogously, for a 3-tensor $S\in\R^{d\times d\times d}$ and a vector $v\in\R^d$, the tensor-vector multiplication is carried out  as follows
\ben*
(S\times_1v)_{ik} &=& \sum_{j=1}^d S_{jik}v_j, \hspace{30pt}
(S\times_2v)_{ik} = \sum_{j=1}^d S_{ijk}v_j,\hspace{30pt}
(S\times_3v)_{ik} = \sum_{j=1}^d S_{ikj}v_j,
\eeqn*
each case resulting in a $d\times d$ matrix. For  example, if $S\in\R^{3\times3\times3}$ and $v\in\R^{3}$ such that 
{\small
\[\begin{tikzpicture}
 \node at (-1.5,-2)      (S)     {$S = $};

\matrix (A) at (0,0) [matrix of math nodes, inner sep=0pt, nodes={inner sep=3pt}, nodes in empty cells,right delimiter={]},left delimiter={[} ]{
S_{111} &  S_{121} & S_{131} \\
S_{211}  & S_{221} & S_{231} \\
S_{311}  & S_{321} & S_{331}  \\
} ;
    \matrix (B) at (1,-1.5) [matrix of math nodes, inner sep=0pt, nodes={inner sep=3pt}, nodes in empty cells,right delimiter={]},left delimiter={[} ]{
S_{112} &  S_{122} & S_{132} \\
S_{212}  & S_{222} & S_{232} \\
S_{312}  & S_{322} & S_{332}  \\
} ;
    \matrix (C) at (2,-3) [matrix of math nodes, inner sep=0pt, nodes={inner sep=3pt}, nodes in empty cells,right delimiter={]},left delimiter={[} ]{
S_{113} &  S_{123} & S_{133} \\
S_{213}  & S_{223} & S_{233} \\
S_{313}  & S_{323} & S_{333}  \\
} ;
    \draw[-] (-1.45,-0.745)-- (0.55,-3.745);
    \draw[-] (1.45,0.735)-- (3.45,-2.265);
\end{tikzpicture}\]
}
then
$$S\times_1 v=\begin{bmatrix} 
S_{111}v_1+S_{211}v_2+S_{311}v_3 & S_{112}v_1+S_{212}v_2+S_{312}v_3 & S_{113}v_1+S_{213}v_2+S_{313}v_3\\
S_{121}v_1+S_{221}v_2+S_{321}v_3 & S_{122}v_1+S_{222}v_2+S_{322}v_3 & S_{123}v_1+S_{223}v_2+S_{323}v_3\\
S_{131}v_1+S_{231}v_2+S_{331}v_3 & S_{132}v_1+S_{232}v_2+S_{332}v_3 & S_{133}v_1+S_{233}v_2+S_{333}v_3
\end{bmatrix}.$$
Generalizing this to arbitrary order tensors yields the following definition.

\begin{definition}[$n$-mode product of a tensor]
The {\boldmath\bf $n$-mode product} of a $k$-order tensor $T\in\R^{d^k}$ with a vector $v\in\R^d$, denoted by $T\times_n v$, is defined elementwise as \[(T\times_n v)_{i_1,\dots,i_{n-1},i_{n+1},\dots,i_k} = \sum_{j=1}^d T_{i_1,\dots,i_{n-1},j,i_{n+1},\dots,i_k} v_j.\]
Note that $T\times_n v\in\R^{d^{k-1}}$, so the order of the resulting tensor is decreased by 1.
\end{definition}

The above definition can also be generalized for tensor-matrix multiplication \cite{tensor1}. Finally we note that the Frobenius norm for matrices can be generalized to tensors in the following way.

\begin{definition}[Tensor Frobenius Norm \cite{tensor1}] The {\bf Frobenius norm of a tensor} $T\in\R^{d^k}$ is the square root of the sum of the squares of all its elements
\[\|T\|_F = \sqrt{\sum_{i_1=1}^d\cdots\sum_{i_k=1}^d {T_{i_1,\dots,i_k}}^2}.\]
\end{definition}

Moments of a distribution have the special property in that they are symmetric in the following sense.

\begin{definition}[Symmetric Tensor]
A tensor $T\in\R^{d^k}$ is {\bf symmetric}, if the tensor is invariant to permutations of the indices, i.e. $$T_{i_1\cdots i_k} = T_{p(i_1\cdots i_k)}$$ for any permutation $p.$
\end{definition}

\begin{remark}
If a tensor is symmetric then the $n$-mode product is independent of the mode, i.e. if $T\in\R^{d^k}$ is symmetric then
\[T\times_n v = T\times_m v\] for any $1\le n,m \le k$.
\end{remark}

The next lemma shows that the tensor Frobenius norm has a particularly simple formula for rank-1 tensors.

\begin{lemma}\label{qual3a} Let $v\in\R^d$ and $k$ be a positive integer then the tensor Frobenius norm of the $k$th-order tensor product is the same as the Euclidean norm of $v$ raised to the $k$, i.e.
\[\|v^{\otimes k}\|_F = {\|v\|}^k.\]
\end{lemma}

\begin{proof}
By the definition of the tensor Frobenius norm, \[
\|v^{\otimes k}\|_F^2 = \sum_{i_1=1}^d\cdots\sum_{i_k=1}^d [(v^{\otimes k})_{i_1,\dots,i_k}]^2 \]
and since $(v^{\otimes k})_{i_1\dots i_k} = v_{i_1}v_{i_2}\cdots v_{i_k}$, we have $\|v^{\otimes k}\|_F^2 = \sum_{i_1=1}^d\cdots \sum_{i_k=1}^d v_{i_1}^2\cdots v_{i_k}^2,$ so
\[\|v^{\otimes k}\|_F^2 = \sum_{i_1=1}^dv_{i_1}^2\sum_{i_2=1}^dv_{i_2}^2\cdots\sum_{i_k=1}^d v_{i_k}^2 = \underbrace{{\|v\|}^2{\|v\|}^2\cdots{\|v\|}^2}_{k \ \rm{times}} \]
by definition of $\|v\|$ so
\[\|v^{\otimes k}\|_F = {\|v\|}^k.\]
\end{proof}

\subsection{Tensor Eigenvectors and Normalized Power Iteration}\label{tensoreigs}

The key to our approximate rank-1 decomposition is rank-1 approximation which is based on tensor eigenvectors which can be found with the Higher Order Power Method (HOPM) which we review in this section.

\begin{definition}[Tensor Eigenvectors and Eigenvalues]
Let $T\in\R^{d^k}$ be a symmetric tensor then $v\in\R^d$ is an {\bf eigenvector} and $\lambda\in\R$ is the corresponding {\bf eigenvalue} of $T$ if \[\left(\left((T\times_1v)\times_1v\right)\cdots \times_1v\right)=\lambda v.\]
\end{definition}
Note that since $T$ is symmetric, the choice of $n$-mode product does not affect the definition of a tensor eigenvector. The next lemma shows that an eigenvalue-eigenvector pair provides a rank-1 approximation of a tensor in the Frobenius norm.
\begin{lemma}\label{eigenapprox}
Let $T$ be a $k$-order symmetric tensor with dimension $d$, i.e. $T\in\R^{d^k}$ and $v\in\R^d$ be a unit eigenvector of $T$ with eigenvalue $\lambda\neq0$. Then   \[\|T-\l v^{\otimes k}\|_F^2 = \|T\|_F^2-\l^2\] and $\|T\|_F\ge\lambda.$
\end{lemma}

\begin{proof}
We first wish to show that $\|T-\l v^{\otimes k}\|_F^2 = \|T\|_F^2-\l^2$.
\ben* 
\|T-\l v^{\otimes k}\|_F^2 &=& \sum_{i_1=1}^d\cdots\sum_{i_k=1}^d [(T-\l v^{\otimes k})_{i_1,\dots,i_k}]^2
= \sum_{i_1=1}^d\cdots\sum_{i_k=1}^d [T_{i_1,\dots,i_k}-\l(v^{\otimes k})_{i_1,\dots,i_k}]^2\\
&=& \sum_{i_1=1}^d\cdots\sum_{i_k=1}^d (T_{i_1,\dots,i_k}^2-2\l T_{i_1,\dots,i_k}v_{i_1}v_{i_2}\cdots v_{i_k}+\l^2 v_{i_1}^2v_{i_2}^2\cdots v_{i_k}^2)\\
&=& \|T\|_F^2-2\l \sum_{i=1}^dv_i(T\times_2v\times_3v\times_4\cdots\times_kv)_i+\l^2\|v^{\otimes k}\|_F
\eeqn*
Since $\|v\|=1$ and by Lemma \ref{qual3a}, $\|v^{\otimes k}\|_F=1$, hence
\ben*
\|T-\l v^{\otimes k}\|_F^2 
&=& \|T\|_F^2-2\l \langle v,\l v\rangle+\l^2 = \|T\|_F^2-2\l^2 \|v\|_2^2+\l^2 
= \|T\|_F^2-\l^2
\eeqn*
Since $\|T-\l v^{\otimes k}\|_F\ge 0$, $\|T\|_F^2-\l^2\ge 0$ so $\|T\|_F^2\ge \l^2$ and taking square roots, 
$\|T\|_F\ge |\l|$.
\end{proof}
It immediately follows from Lemma \ref{eigenapprox} that the eigenvector with the largest eigenvalue will achieve the best rank-1 approximation among the eigenpairs. In fact, it has been shown that the eigenpair with the largest eigenvalue achieves the best possible rank-1 approximation of the tensor \cite{EigApprox,HOPM}. This fact will form the basis for an effective algorithm for finding an approximate rank-1 decomposition in the next section.

Finally, an effective algorithm for finding the eigenvector associated to the largest eigenvalue in absolute value is the Higher Order Power Method (HOPM) originally developed in \cite{HOPMoriginal} and further analyzed in \cite{HOPMrevisited,HOPM}.  In the case of symmetric tensors the Symmetric-HOPM (S-HOPM) has a simpler form that is very similar to Normalized Power Iteration (NPI) but is not guaranteed to converge \cite{EigApprox}.  The HOPM algorithm for a symmetric order-k tensor $T \in \R^{d^k}$ requires initialization with the left singular vector, $u$, corresponding to the largest singular value of the unfolding (reshaping) of the tensor into a $d\times d^{k-1}$ matrix.  The HOPM then defines a sequence of vectors inductively by, $v_0^{(1)}=\cdots =v_0^{(k)} = u$ and then sequentially updates 
\begin{align} w &= T \times_1 v_{j+1}^{(1)} \times_1 \cdots \times_1 v_{j+1}^{(i-1)} \times_1 v_j^{(i+1)} \times_1 \cdots \times_1 v_j^{(k)} \nonumber \\
v_{j+1}^{(i)} &= \frac{w}{||w||} \nonumber
\end{align}
for each $i=1,...,k$ and then increments $j$.
Notice that the product that updates $v_{j+1}^{(i)}$ is the tensor $T$ multiplied by the $k-1$ other vectors, leaving out $v_j^{(i)}$.  Also note that we use the already updated $(j+1)$-step vectors for the first $i-1$ products and the $j$-step vectors for the last $k-i$ products.  The HOPM is guaranteed to converge to an eigenvector of $T$ \cite{HOPMrevisited}, and when $T$ is symmetric all $v^{(1)}_j,...,v^{(k)}_j$ converge to the same eigenvector but may differ in sign for even order tensors.  

For completeness we summarize the HOPM algorithm of \cite{HOPMoriginal}.
\begin{algorithm}[H] 
\caption{Higher Order Power Method (HOPM) \cite{HOPMoriginal} } 
\label{algHOPM} 
\begin{algorithmic} 
	\State {\bf Inputs:} A $k$-tensor $T \in \mathbb{R}^{d^k}$
	\State {\bf Outputs:} Eigenvector $v\in\mathbb{R}^d$ and eigenvalue $\lambda$ such that $T\times_1 v \times_1 \cdots \times_1 v = \lambda v$
	\State
	\Indent
	\State Reshape $T$ into a $d \times d^{k-1}$ matrix and compute the leading left singular vector, $v_0$
	\State Initialize $v^1=v^2=\cdots=v^k=v_0$, $\lambda = {\rm Inf}$ and $\lambda_{\rm prev}=0$
	\While{$|\lambda-\lambda_{\rm prev}| > {\rm tol}$}
	    \For{$\ell=1,...,k$}
	        \State Set $v^\ell_s = \sum_{i_1,...,i_{\ell-1},i_{\ell+1},...,i_k=1}^d T_{i_1,...,i_{\ell-1},s,i_{\ell+1},...,i_k}v_{i_1}^1 \cdots v_{i_{\ell-1}}^{\ell-1} v_{i_{\ell+1}}^{\ell+1} \cdots v_{i_k}^k$
	    \EndFor
	    \State Set $\lambda_{\rm prev} = \lambda$
	    \State Set $\lambda = \sum_{i_1,...,i_k=1}^d T_{i_1,...,i_k}v^1_{i_1}\cdots v^1_{i_k}$
	\EndWhile
	\State Set $v=v^1$
	\EndIndent
	\State Return $v, \lambda$.
    \State 
\end{algorithmic}
\end{algorithm}

Unlike the case of matrices, for tensors of order greater than two the basins of attraction for multiple distinct eigenvalues can have non-zero measure.  It has been observed \cite{HOPM,HOPMrevisited,EigApprox} that initialization with the left singular vector, $u$, of the tensor unfolding typically leads to convergence to the eigenvector with the largest eigenvalue.  The next section will rely on the ability to find the eigenpair associated to the largest eigenvalue (in absolute value) so a guaranteed way to find an initial condition in the basin of the largest eigenvalue is still an important problem for future research.

\section{Approximate Rank-1 Decomposition}\label{approxrank1}


In this section we show how tensors eigenvectors can be used to form an approximate rank-1 decomposition up to an arbitrary level of precision. Of course, this is not a method of finding the minimal rank-1 decomposition, the computation of which is NP-complete \cite{NP1,NP2}. Moreover, we do not even see an exact rank-1 decomposition.  Instead, given an order-$k$ tensor $T$, we seek a sequence of vectors $v_\ell$ and constants $\lambda_\ell$ such that $\sum_{\ell=1}^L \lambda_\ell v_\ell^{\otimes k}$ approximates $T$ in the Frobenius norm up to an error that can be made arbitrarily small by increasing $L$.  In the next section we will show that this approximate rank-1 decomposition is a key component for generalizing the unscented ensemble to higher moments.

Our approach is motivated by a theorem of \cite{EigApprox} which states that if $v$ is the unit eigenvector of an order-$k$ tensor $T$ associated to the largest eigenvalue $\lambda$ (in absolute value), then $\lambda v^{\otimes k}$ is the best rank-1 approximation of $T$, namely \[\|T-\lambda v^{\otimes k}\|\] is minimized over all possible $\lambda$, $\|v\|=1$. It is well known that subtracting the best rank-1 approximation does not produce an \emph{exact} rank-1 decomposition, and in fact may increase tensor rank \cite{subtract1,subtract2}.  However, it was suggested in \cite{EigApprox} that repeatedly subtracting the rank-1 approximations may result in an \emph{approximate} rank-1 decomposition. The following theorem shows that this process converges subject to a certain tensor eigenvalue inequality that will be shown in Lemma \ref{order3and4} below.
	

\vspace{1em}
\begin{theorem}\label{rank1decomp}
Let $T$ be a $k$-order symmetric tensor with size $d$, i.e. $T\in\R^{d^k}$. Consider the process of finding an approximate rank-1 decomposition of $T$ by starting from $T_0 = T$ and setting $T_{\ell+1} = T_{\ell}-\lambda_\ell v_\ell^{\otimes k}$ where $\lambda_\ell$ is the largest eigenvalue in absolute value of $T_\ell$ and $v_\ell$ is the associated eigenvector. Assume also that there exists a universal constant $c\in(0,1]$ such that $\l_{\ell}\ge c|{(T_{\ell})}_{i_1\dots i_k}|$. Then $\|T_\ell\|_F \to 0$ and for $\displaystyle r = \sqrt{1-\frac{{c}^2}{d^k}}\in(0,1)$ 
\[\frac{\|T_{\ell+1}\|_F}{\|T_\ell\|_F}\le r \hspace{1em} \hbox{and}\hspace{1em} T = \sum_{\ell=1}^L\l_\ell v_\ell^{\otimes k}+\mathcal{O}(r^L)\]
for all $L\in\N$.
\end{theorem}

\begin{proof}
First let $\lambda_{maxabs}$ be the largest eigenvalue in absolute value of a tensor $T$ and assume $\l_{maxabs}\ge c|T_{i_1\dots i_k}|$ for all $i_1,\dots,i_k$. We will show that there exists a constant $c_2=\frac{c}{d^{k/2}}\in(0,1]$ such that $\l_{maxabs}\ge c_2\|T\|_F$. Since $\l_{maxabs}\ge c|T_{i_1\dots i_k}|$, we have \[\l_{maxabs}^2\ge c^2T_{i_1\dots i_k}^2\] which implies that \[d^k\l_{maxabs}^2\ge c^2\sum_{i_1,\dots,i_k}T_{i_1\dots i_k}^2\] so we have $d^{k/2}\l_{maxabs}\ge c\sqrt{\sum_{i_1,\dots,i_k}T_{i_1\dots i_k}^2}$ and  
\begin{eqnarray}\l_{maxabs}\ge \frac{c}{d^{k/2}}\|T\|_F,\label{c2ineq}\end{eqnarray} where we take $\displaystyle c_2=\frac{c}{d^{k/2}}\in(0,1)$, since $c\in(0,1)$ and $d\ge1$. 
By Lemma \ref{eigenapprox} applied to $T_\ell$, we have \[{\|T_{\ell+1}\|_F}^2={\|T_\ell-\l_\ell v_\ell^{\otimes k}\|_F}^2 = {\|T_\ell\|_F}^2-{\l_\ell}^2.\] Since $\lambda_\ell$ is defined to be the largest eigenvalue of $T_\ell$  \eqref{c2ineq} says that $\l_\ell\ge c_2\|T_\ell\|_F$ where $c_2=\frac{c}{d^{k/2}}$ so 
\ben*
{\|T_{\ell+1}\|_F}^2 &\le&{\|T_\ell\|_F}^2-{c_2}^2{\|T_\ell\|_F}^2\\ &\le&(1-{c_2}^2){\|T_\ell\|_F}^2.
\eeqn*
Thus, setting $r=\sqrt{1-{c_2}^2}\in(0,1)$ we have $\|T_{\ell+1}\|_F \le r\|T_\ell\|_F$ and $\|T_{\ell+1}\|_F \le r^2\|T_{\ell-1}\|_F$ and so forth and proceeding inductively we find, \[\|T_{\ell+1}\|_F \le r^{\ell+1}\|T_{0}\|_F = r^{\ell+1}\|T\|_F.\] Since $0<r<1$, $\lim_{\ell\to\infty} r^{\ell+1}=0$, so $0\le \|T_{\ell+1}\|_F\le r^{\ell+1}\|T\|_F \to 0$ implies $\|T_{\ell+1}\|\to0$ as $\ell\to\infty$. Since this limit is 0, an upper bound on the rate of convergence of $\|T_\ell\|_F$ is found by considering \[\frac{\|T_{\ell+1}\|_F}{\|T_\ell\|_F}\le r = \sqrt{1-\frac{{c}^2}{d^k}}.\]

\end{proof}
Theorem \ref{rank1decomp} gives an effective algorithm for finding approximate rank-1 decompositions of tensors, however it requires an inequality of the form 
\begin{align}\label{eigineq} \lambda_{maxabs}\ge c|T_{i_1,\dots,i_k}|.\end{align}
The inequality \eqref{eigineq} holds for symmetric matrices with $c=1$, since if $T \in \mathbb{R}^{d^2}$ is symmetric it has an orthogonal eigendecomposition $T = U^\top \Lambda U$ so by the Cauchy-Schwarz inequality, 
\[ |T_{ij}| = \left|\left<u_i,\lambda_j u_j\right>\right| \leq ||u_i||\,||\lambda_j u_j|| = |\lambda_j| \leq \lambda_{maxabs} \]
and the identity matrix shows that $c=1$ is the best possible constant for matrices.  Of course, this method of proof cannot be generalized to arbitrary tensors due to the lack of a similar rank-1 eigendecomposition.  Nevertheless, the next lemma shows that an inequality of the form \eqref{eigineq} does hold for all symmetric tensors of orders 3 and 4.
\begin{lemma}\label{order3and4}
If $T$ is a symmetric 3-tensor, then \[\lambda_{maxabs}\ge \frac{2}{3 + 4 \sqrt{2} + \sqrt{3}}|T_{ijk}|.\] If $T$ is a symmetric 4-tensor, then \[\lambda_{maxabs}\ge \frac{6}{323}|T_{ijk\ell}|.\]
\end{lemma}
The proof of Lemma \ref{order3and4} is quite involved and can be found in the appendix. We conjecture that such an inequality holds for symmetric tensors of any order, and we note that the constants in Lemma \ref{order3and4} are not known to be sharp.  For the purposes of this paper, we are focused on matching the skewness and kurtosis of a distribution so we only need the approximate rank-1 decomposition for tensors up to order 4. In the next section we will show how to use the approximate rank-1 decomposition to build an ensemble that simultaneously matches the mean, covariance, skewness and kurtosis.

\begin{figure}[H]
  \centering
  \subfloat[]{\includegraphics[width=0.245\textwidth]{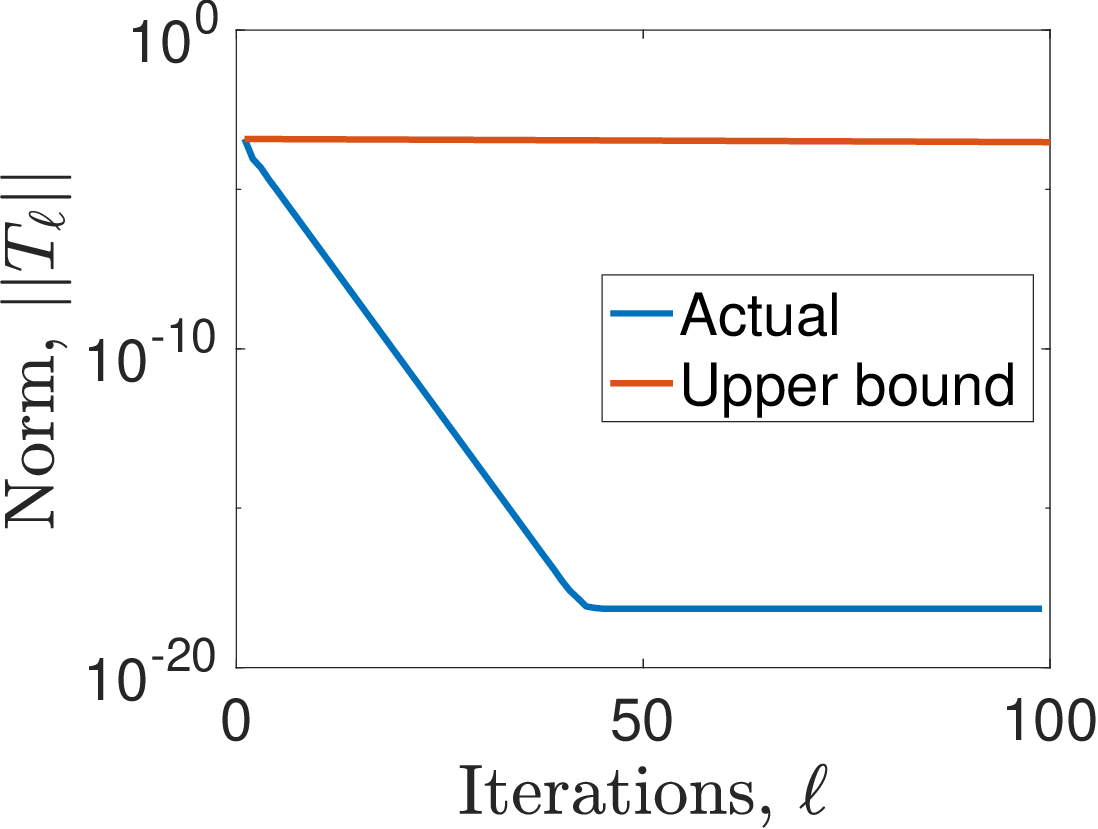}\label{fig:f1a}}
  \hfill
  \subfloat[]{\includegraphics[width=0.245\textwidth]{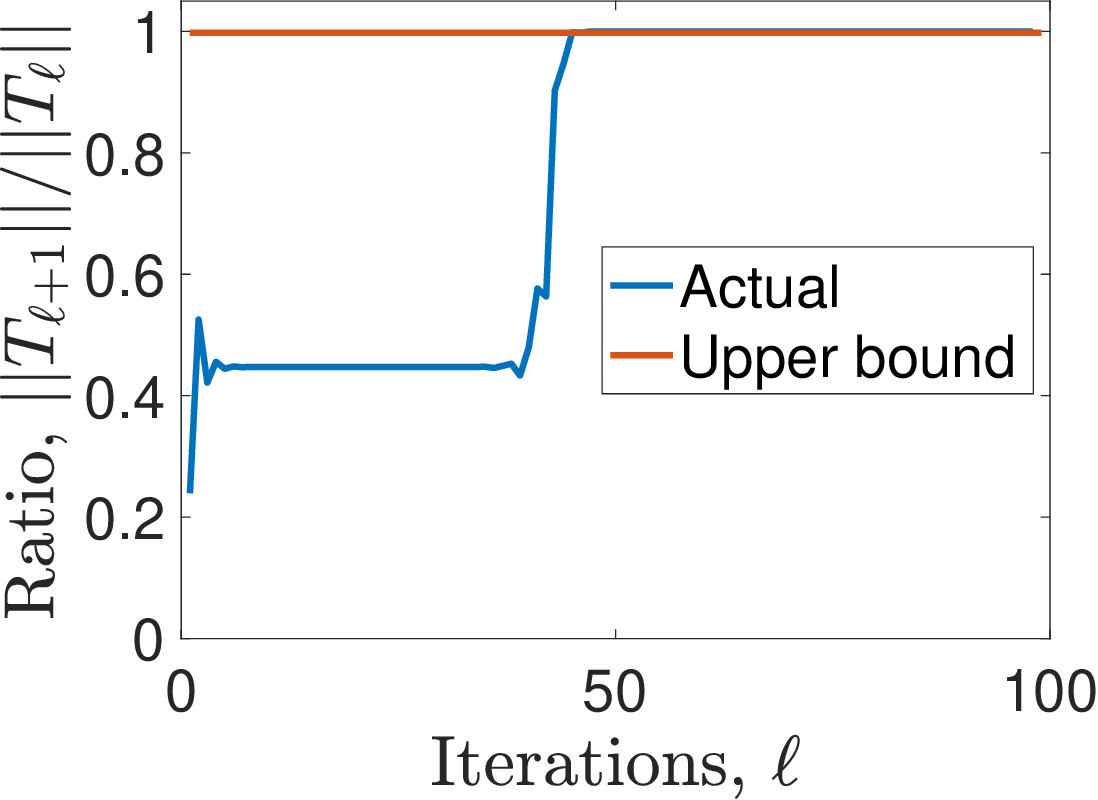}\label{fig:f1b}}
  \hfill
  \subfloat[]{\includegraphics[width=0.245\textwidth]{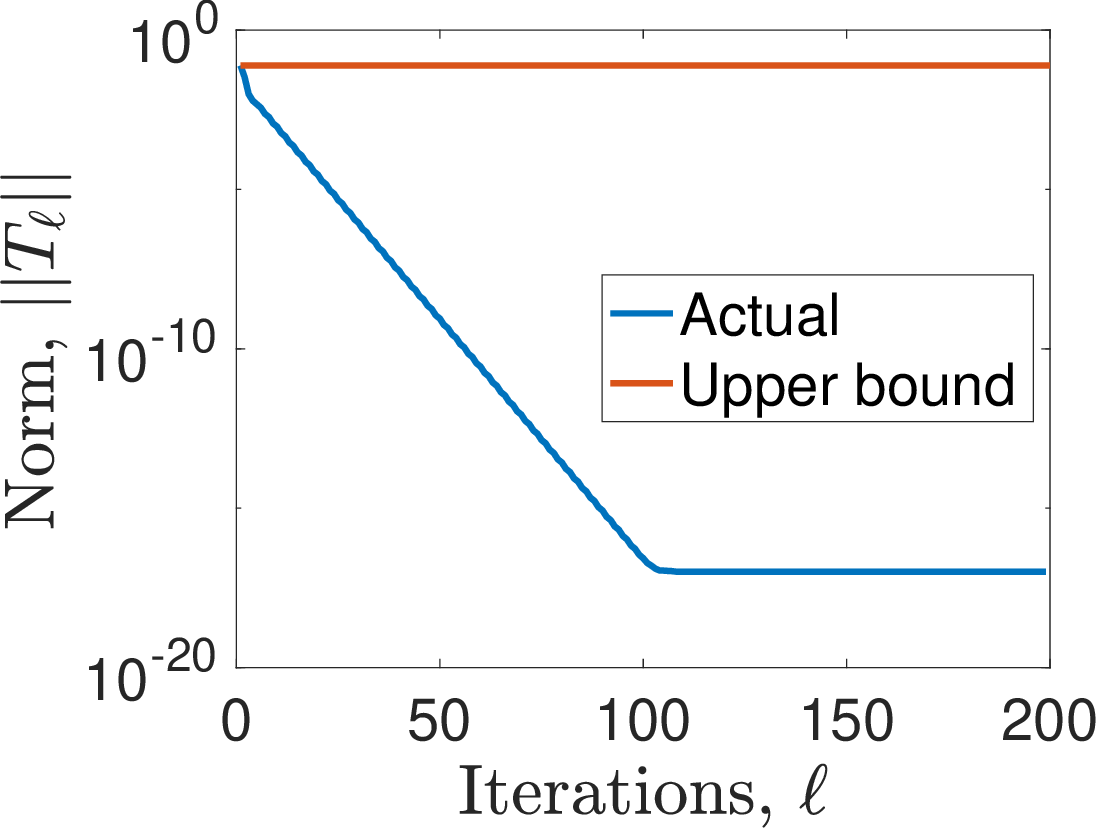}\label{fig:f1c}}
  \hfill
  \subfloat[]{\includegraphics[width=0.245\textwidth]{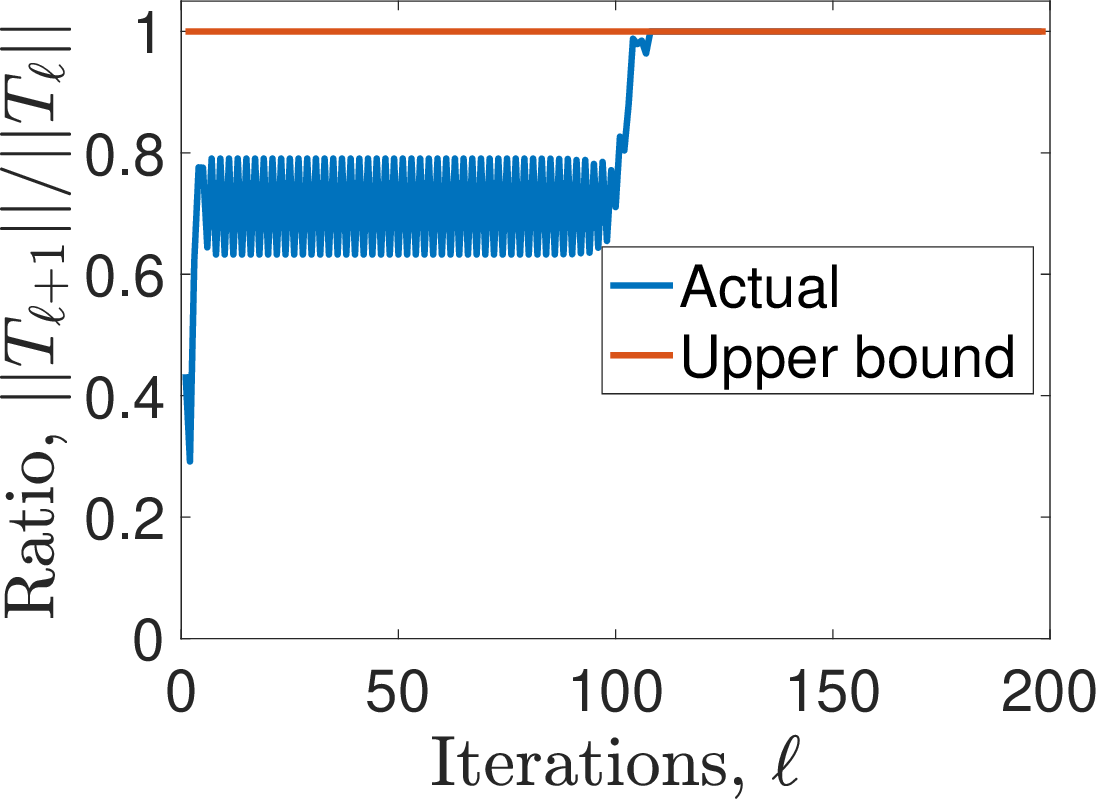}\label{fig:f1d}}
  
  \subfloat[]{\includegraphics[width=0.245\textwidth]{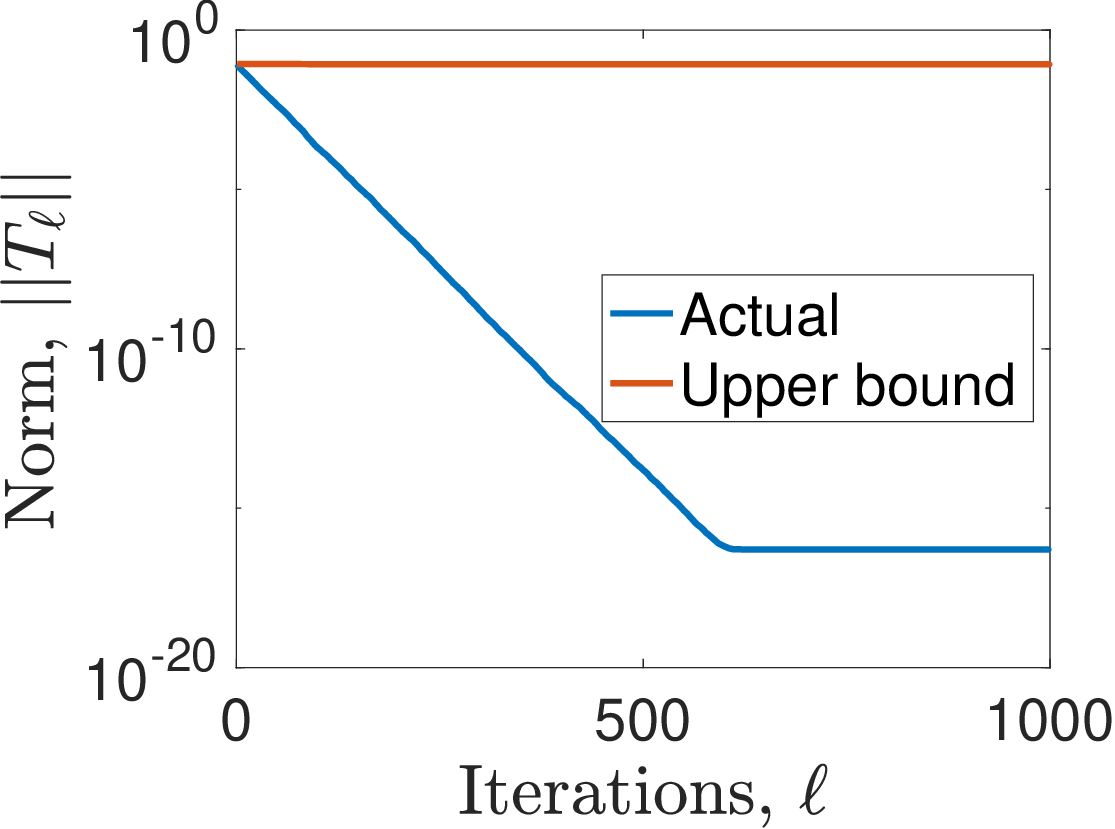}\label{fig:f1e}}
  \hfill
  \subfloat[]{\includegraphics[width=0.245\textwidth]{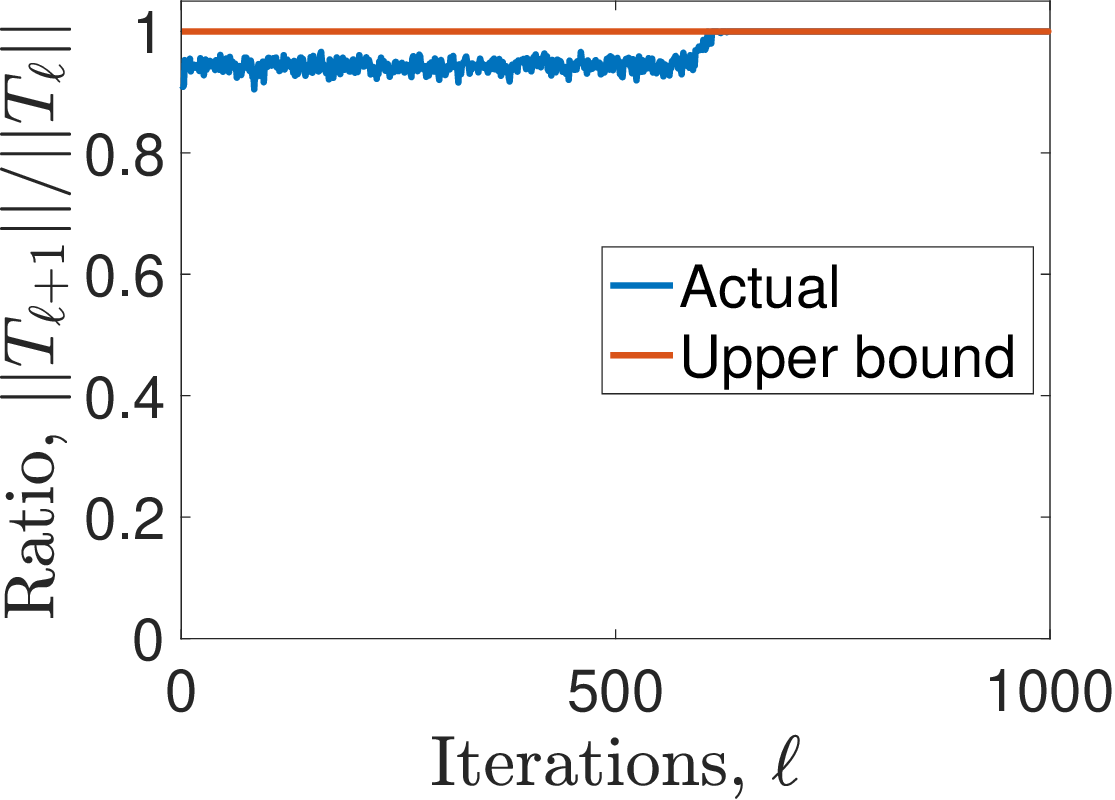}\label{fig:f1f}}
  \hfill
  \subfloat[]{\includegraphics[width=0.245\textwidth]{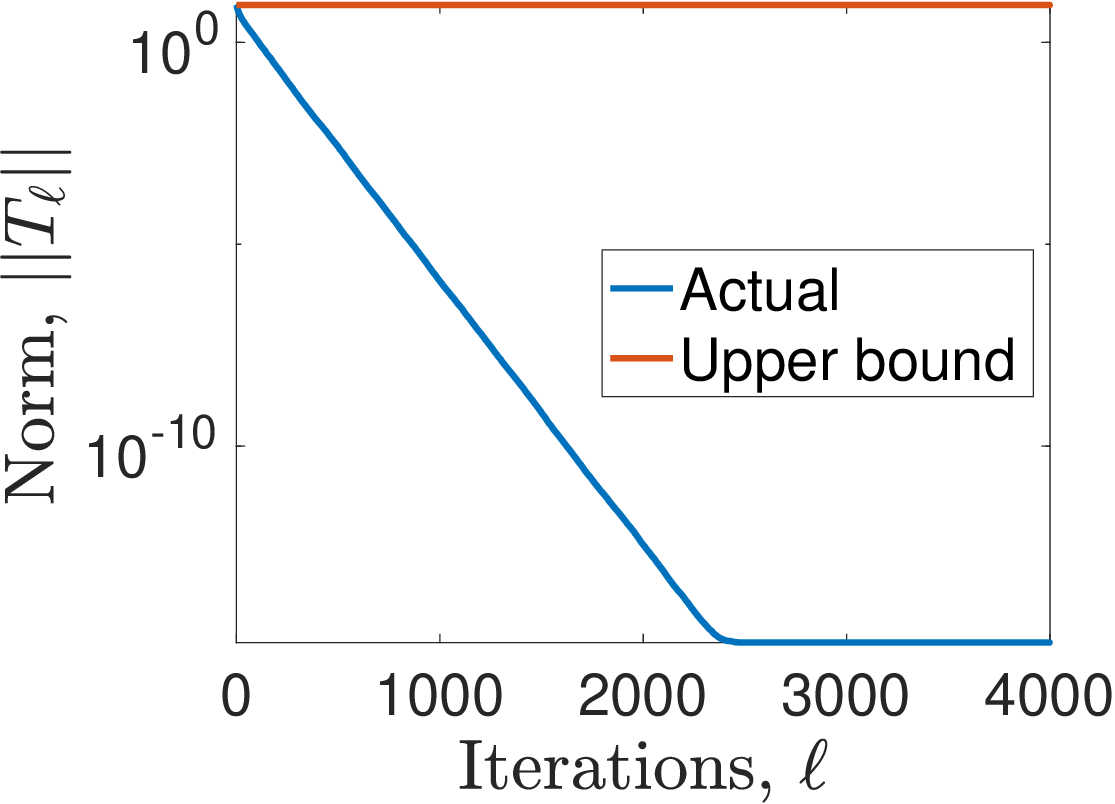}\label{fig:f1g}}
  \hfill
  \subfloat[]{\includegraphics[width=0.245\textwidth]{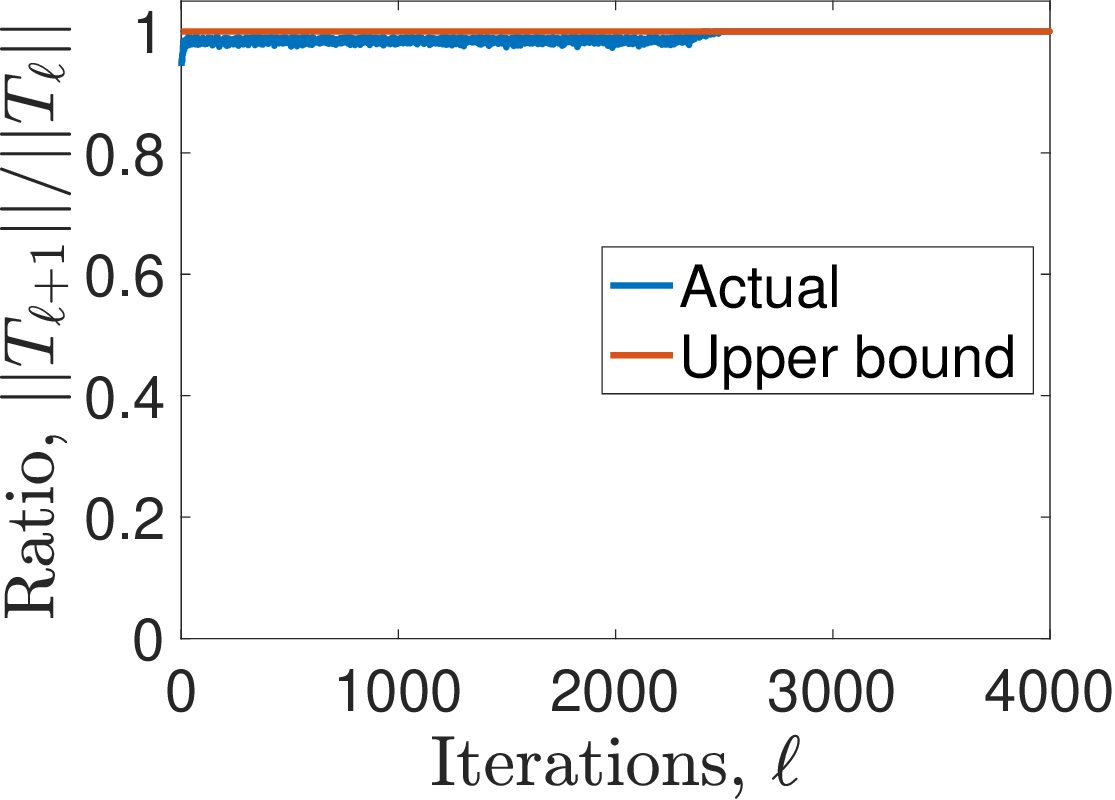}\label{fig:f1h}}
  \caption{Top (a-d): With $d=2$ we demonstrate the approximate rank-1 decomposition on a random 3-tensor (a,b) and 4-tensor (c,d).  The norm of the residual (a,c)(blue) decays to numerical zero faster than the theoretical upper bound, $r^{\ell}$ (red). The ratio of successive Frobenius norms (b,d)(blue) is always less than the derived upper bound, $r$ (red). Bottom (e-h): We repeat the experiment with $d=10$.}\label{fig1}
\end{figure}

We summarize the approximate rank-1 decomposition algorithm below.
\begin{algorithm}[H] 
\caption{Approximate Rank-1 Decomposition } 
\label{algRank1} 
\begin{algorithmic} 
	\State {\bf Inputs:} A $k$-tensor $T \in \mathbb{R}^{d^k}$ and a tolerance $\tau$.
	\State {\bf Outputs:} Vectors, $v_i$, and signs, $s_i \in \{-1,1\}$ such that $\left|\left| \sum_{i=1}^J s_i v_i^{\otimes k} - T \right|\right|_F \leq \tau$.
	\State
	\Indent
	\State Set $i=1$
	\While{$||T||_F > \tau$}
	    \State Apply the HOPM (Algorithm \ref{algHOPM}) to find an eigenpair $(v,\lambda)$ of $T$.
	    \State Set $s_i = {\rm sign}(\lambda)$ (note that if $k$ is odd we can always choose $s_i = 1$)
	    \State Set $v_i = |\lambda|^{1/k}v$
	    \State Set $T = T - s_i v_i^{\otimes k}$
	    \State Set $i = i+1$
	\EndWhile
	\EndIndent
	\State Return the set of all $s_i, v_i$.
    \State 
\end{algorithmic}
\end{algorithm}
Finally, we demonstrate this algorithm on a random 3-tensor and 4-tensor with $d=2$ and $d=10$ in Figure \ref{fig1}. We note that in all cases the convergence is much faster than our theoretical upper bound, however for $d=10$ we see that the ratio of residual norms approaches much closer to our upper bound.  Moreover, high dimensional tensors require a much larger number of vectors to achieve a given tolerance with the approximate rank-1 decomposition.  So while our approach provides an effective solution, it is likely that there is room for improvement, and the Higher Order Unscented Transform (HOUT) introduced in the next section can use any method of rank-1 decomposition. 

\section{Higher Order Unscented Transform}\label{sec:HOUT} The goal of the scaled unscented transform is to generate a small ensemble that exactly matches the mean and covariance of a distribution, thus forming a quadrature rule that can be to estimate the expected value of nonlinear functions.  In this section we define the Higher Order Unscented Transform which matches the first four moments of a distribution, thus providing a quadrature rule with a higher degree of exactness.  While we only describe the process explicitly for up to four moments, our method is based on the approximate tensor decomposition from the previous section and should allow generalization to an arbitrary number of moments.

Suppose we are given the following moments of the distribution of a random variable: the mean $\mu\in\R^d$, the covariance matrix $C\in\R^{d\times d}$, the skewness tensor $S\in\R^{d\times d\times d}$, and kurtosis tensor $K\in\R^{d\times d\times d\times d}$.  Let $\tau$ be a parameter that specifies the tolerance of the approximate rank-1 decompositions and let $S$ and $K$ have the approximate rank-1 decompositions 
\[ \left|\left|S - \sum\limits_{i=1}^J \tilde{v_i}^{\otimes 3}\right|\right|_F \leq \tau/2 \hspace{30pt} \left|\left| K - \sum\limits_{i=1}^L s_i\tilde{u_i}^{\otimes 4}\right|\right|_F \leq \tau/2 \]
where $s_i \in \{-1,1\}$ denote signs.  Note that these approximate decompositions can be constructed by the algorithm described in Theorem \ref{rank1decomp} and then moving the eigenvalues inside the tensor power by the rule $(cv)^{\otimes k} = c^k v^{\otimes k}$.  Note that the signs $s_i$ are required for the kurtosis since constants come out of even order tensor powers as absolute values.

The key to forming an ensemble that matches all four moments simultaneously is carefully balancing the interactions between the moments.  For example, if we add new quadrature nodes of the form $\mu+\gamma\tilde v_i$ in order to try to match the skewness, these nodes will influence the mean of the ensemble.  In order to balance these interactions we make the following definitions based on the approximate rank-1 decompositions of the skewness and kurtosis,
\[ \tilde{\mu} = \sum\limits_{i=1}^J \tilde{v_i}, \hspace{40pt} \displaystyle\hat{\mu} = - \gamma^{-2}\tilde{\mu}, \hspace{40pt} \tilde{C} = \sum\limits_{i=1}^L s_i\tilde{u_i}^{\otimes 2}, \hspace{40pt} \displaystyle\hat{C}=C-\frac1{\delta^2}\tilde{C} \] 
where $\hat{L}=\sum\limits_{i=1}^Ls_i$ and $\beta,\gamma,\delta$ are arbitrary positive constants that will define the 4 moment $\sigma$-points below.  We note that $C$ is assumed symmetric and positive definite since it is a covariance matrix and $\tilde C$ is symmetric by definition.  In order to insure that $\hat C$ is also positive definite, let $\lambda_{\max}^{\tilde C}$ be the largest eigenvalue of $\tilde C$ and let $\lambda_{\rm min}^C$ be the smallest eigenvalue of $C$, then we require that $\delta>\sqrt{\frac{\lambda_{\rm max}^{\tilde{C}}}{\lambda_{\rm min}^C}}$ which guarantees that $\hat C$ is positive definite.  We note that this choice can be overly conservative especially when $C$ is close to rank deficient.  In these cases, it can be helpful to iterative divide $\delta$ by 2 as long as $\hat C$ remains positive definite.  These choices balance out the interactions between the moments and are the key to proving Theorem \ref{fourmoments} below.  We are now ready to define the 4-moment $\sigma$-points.

\begin{definition}[The 4-moment $\sigma$--points of the Higher Order Unscented Transform]\label{sigmapts}

Let $\alpha$, $\b$, $\gamma$, $\delta$ be positive real numbers, we define the {\boldmath{\bf$4$ moment $\sigma$-points}} by 
\[\sigma_i = 
\begin{cases} 
    \mu & \hbox{if } \  i=-2 \\
    \mu+\alpha\hat{\mu} & \hbox{if } \  i=-1 \\
    \mu-\alpha\hat{\mu} & \hbox{if } \  i=0 \\
    \mu+\beta\sqrt{\hat{C}}_i & \hbox{if } \  i=1,\dots, d \\
    \mu-\beta\sqrt{\hat{C}}_{i-d} & \hbox{if } \  i=d+1,\dots, 2d \\
    \mu+\gamma\tilde{v}_{i-2d} & \hbox{if } \ i=2d+1,\dots, 2d+J\\
     \mu-\gamma\tilde{v}_{i-2d-J} & \hbox{if } \ i=2d+J+1,\dots, 2d+2J\\
     \mu+\delta\tilde{u}_{i-2d-2J} & \hbox{if } \ i=2d+2J+1,\dots, 2d+2J+L\\
     \mu-\delta\tilde{u}_{i-2d-2J-L} & \hbox{if } \ i=2d+2J+L+1,\dots, N
   \end{cases}
\]
and the corresponding {\boldmath {\bf weights}} by
\[w_i = 
\begin{cases} 
    1-d\beta^{-2}-\hat L \delta^{-4} & \hbox{if } \  i=-2 \\
    \frac{1}{2}\alpha^{-1} & \hbox{if } \  i=-1 \\
    -\frac{1}{2}\alpha^{-1} & \hbox{if } \  i=0 \\
    \frac{1}{2}\beta^{-2} & \hbox{if } \  i=1,\dots, 2d \\
    \frac{1}{2}\gamma^{-3} & \hbox{if } \ i=2d+1,\dots, 2d+J\\
    -\frac{1}{2}\gamma^{-3} & \hbox{if } \ i=2d+J+1,\dots, 2d+2J\\
    \frac{1}{2}\delta^{-4}s_{i-2d-2J} & \hbox{if } \ i=2d+2J+1,\dots, 2d+2J+L\\
    \frac{1}{2}\delta^{-4}s_{i-2d-2J-L} & \hbox{if } \ i=2d+2J+L+1,\dots, N
   \end{cases}
\]
For convenience, denote $N=2(d+J+L)$.
\end{definition}
The next theorem shows that the 4-moment $\sigma$-points match the first two moments exactly and match the skewness and kurtosis up to an error term that can be controlled below. 
\begin{theorem}\label{fourmoments}
Given the 4-moment $\sigma$-points associated with $\mu$, $C$, $S$, and $K$  we have $\sum_{i=-2}^N w_i = 1$ and
\ben*
\sum_{i=-2}^{N} w_i\sigma_i &=& \mu \\
\sum_{i=-2}^{N} w_i(\sigma_i-\mu)^{\otimes 2}  &=& C\\
\left|\left|\sum_{i=-2}^{N} w_i(\sigma_i-\mu)^{\otimes 3} - S \right|\right|_F &\leq& \tau/2 +  \alpha^2\left|\left|\hat{\mu}^{\otimes 3}\right|\right|_F \\
\left|\left|\sum_{i=-2}^{N} w_i(\sigma_i-\mu)^{\otimes 4} - K\right|\right|_F &\leq& \tau/2 + \beta^2 \left|\left|\bar C\right|\right|_F.
\eeqn*
where $\bar C = \sum_{i=1}^{d} \sqrt{\hat{C}}_i^{\otimes 4}$.
\end{theorem}
Notice that the third and fourth moment equations do not exactly match the skewness and kurtosis, respectively. Of course, we only used an approximate rank-1 decomposition to begin with, which accounts for the $\tau$ term in the error. Thus, the real goal is to bound the other error term by the same tolerance, $\tau$. The following corollary shows how to control the error terms on the skewness and kurtosis. 
\begin{corollary}\label{tolerance} Let $\tau$ be a specified tolerance for the absolute error of the skewness and kurtosis and set $\bar{C}=\sum\limits_{i=1}^{d}\sqrt{\hat{C}}_i^{\otimes 4}$ and $\hat \mu$ as in Theorem \ref{fourmoments}. If we choose parameters $\alpha,\beta$ such that 
\[ \alpha < \sqrt{ \frac{\tau}{2||\hat{\mu}^{\otimes 3}||_F}} \hspace{30pt}{\rm and }\hspace{30pt} \beta < \sqrt{ \frac{\tau}{2||\overline{C}||_F }} \] 
then 
\[  \left|\left|\sum_{i=-2}^{N} w_i(\sigma_i-\mu)^{\otimes 3} - S \right|\right|_F < \tau \hspace{30pt}{\rm and }\hspace{30pt} 
\left|\left|\sum_{i=-2}^{N} w_i(\sigma_i-\mu)^{\otimes 4} - K \right|\right|_F  < \tau. \]
\end{corollary}
\begin{proof}  The inequality for $\beta$ follows immediately from Theorem \ref{fourmoments}.  Once $\beta$ is chosen, then we can define,
\begin{align} ||\hat{\mu}^{\otimes 3}||_F &= \left|\left| \left(\left(1-d\beta^{-2}-\hat L\delta^{-4}\right)\mu - \gamma^{-2}\tilde{\mu}\right)^{\otimes 3}\right|\right|_F \nonumber \end{align}
and choosing $\alpha < \sqrt{ \frac{\tau}{2||\hat{\mu}^{\otimes 3}||_F}}$  we have
\[ \left|\left|\sum_{i=-2}^{N} w_i(\sigma_i-\mu)^{\otimes 3} - S \right|\right|_F \leq \tau/2 + \alpha^2|| \hat{\mu}^{\otimes 3}||_F  <\tau\]
as desired.  
\end{proof}
Corollary \ref{tolerance} could easily be reformulated to control relative error if desired, and taken to the extreme we could make the quadrature rule exact up to numerical precision.  As a practical matter, this is not an effective strategy since it would result in a larger condition number for the numerical quadrature as shown in the following remark.

\begin{algorithm}[H]
\caption{\label{HOUTalg} Higher Order Unscented Transform (HOUT) }
\begin{algorithmic} 
	\State {\bf Inputs:} A function $f$, tolerance $\tau$, and the mean, $\mu$, covariance, $C$, skewness, $S$, and kurtosis, $K$, of a random variable $X$.
	\State {\bf Outputs:} Estimate of $\mathbb{E}[f(X)]$ with degree of exactness 4.
	\State
	\Indent
    \State Compute the approximate rank-1 decomposition $\left|\left| S - \sum_{i=1}^J \tilde v^{\otimes 3}\right|\right|_F \leq \tau/2$
    \State Compute the approximate rank-1 decomposition $\left|\left| K - \sum_{i=1}^L  s_i \tilde u^{\otimes 4}\right|\right|_F \leq \tau/2$
    \State Set $\tilde C = \sum_{i=1}^L s_i \tilde u_i^{\otimes 2}$.
    \State Compute the largest eigenvalue $\lambda_{\max}^{\tilde C}$ of $\tilde C$ and the smallest eigenvalue $\lambda_{\min}^{C}$ of $C$
    \State Choose $\delta > \sqrt{\frac{\lambda_{\max}^{\tilde C}}{\lambda_{\min}^{C}}}$ (note that $C$ is positive definite so $\lambda_{\min}^{C}>0$)
    \State (Optional) While $C - \delta^{-2}\tilde C$ is positive definite, set $\delta = \delta/2$
    \State Set $\hat C = C - \delta^{-2}\tilde C$  
    \State Compute the symmetric square root of $\hat C$ with columns $\sqrt{\hat C_i}$
    \State Set $\bar C = \sum_{i=1}^d \sqrt{\hat C_i}^{\otimes 4}$
    \State Choose $\beta < \sqrt{\frac{\tau}{2||\bar C||_F}}$ and choose $\gamma > 0$ (default $\gamma=J^{-1/3}$)
    \State Set $\hat L = \sum_{i=1}^L s_i$ and $\tilde \mu = \sum_{i=1}^J \tilde v_i$ and $\hat \mu = (1-d\beta^{-2}-\hat L \delta^{-4})\mu - \gamma^{-2}\tilde \mu$
    \State Choose $\alpha < \sqrt{\frac{\tau}{2||\hat \mu^{\otimes 3}||_F}}$
    \State Define the 4-moment $\sigma$-points, $\sigma_i$, and weights, $w_i$, according to Definition \ref{sigmapts}
    \EndIndent
    \State Output: $\sum_{i=-2}^N w_i f(\sigma_i)$
\end{algorithmic}
\end{algorithm}

\begin{remark}
The absolute condition number of the Higher Order Unscented Transform is bounded above $\sum_{i=0}^N |w_i|$.  Using the bounds from Corollary \ref{tolerance} we find
\begin{align} \sum_{i=0}^{N}|w_i| &= \frac{1}{\alpha} + \frac{d}{\beta^2} + \frac{J}{\gamma^3} + \frac{L}{\delta^4} > \sqrt{\frac{||\bar{\mu}^{\otimes 3}||_F}{\tau}} + \frac{d ||\bar C||_F}{\tau} +\frac{J}{\gamma^3} + \frac{L}{\delta^4} = \mathcal{O}(\tau^{-1}) \nonumber  \end{align}
which shows that the condition number has the potential to blow up as the tolerance is decreased.
\end{remark}
We summarize the HOUT algorithm in Algorithm \ref{HOUTalg} and we now turn to some numerical experiments to demonstrate the HOUT.

\section{Numerical Experiments}\label{sec:numerics}

\begin{figure}[H]
  \centering
   \subfloat[]{\includegraphics[width=0.35\textwidth]{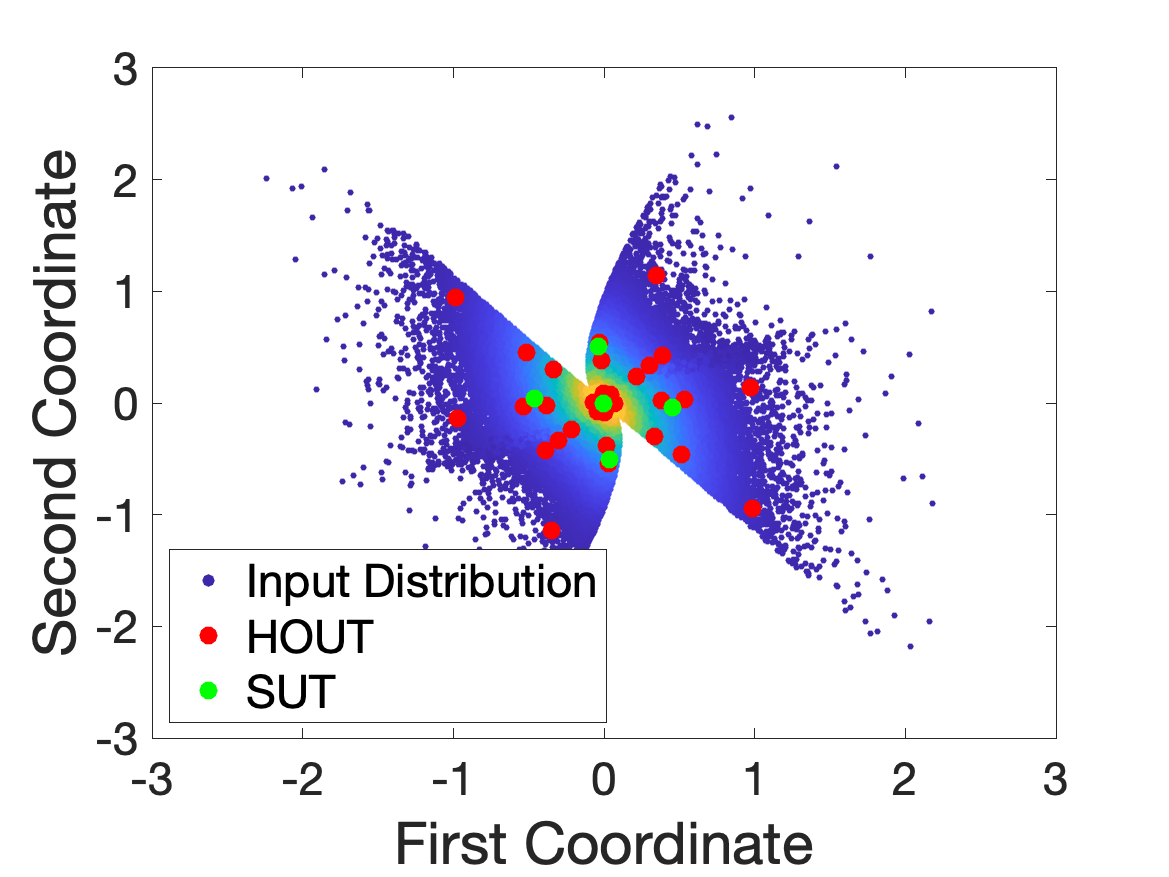}} \hfill \subfloat[]{\includegraphics[width=0.315\textwidth]{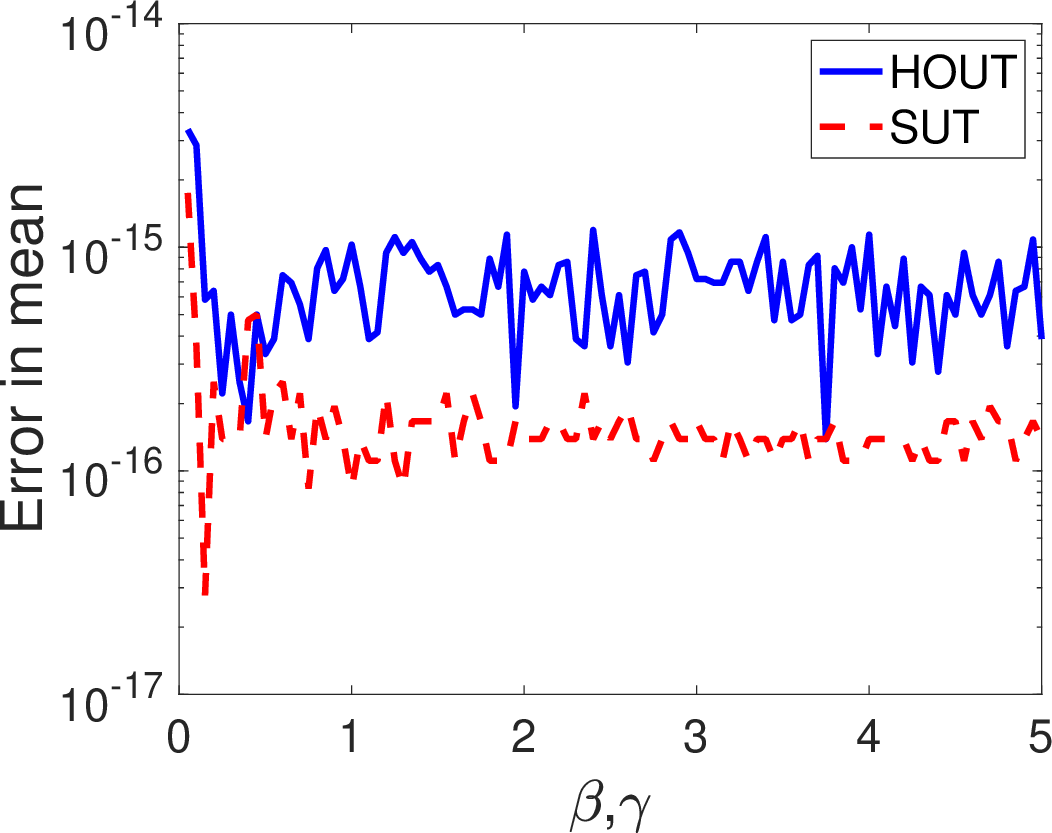}} \hfill \subfloat[]{\includegraphics[width=0.315\textwidth]{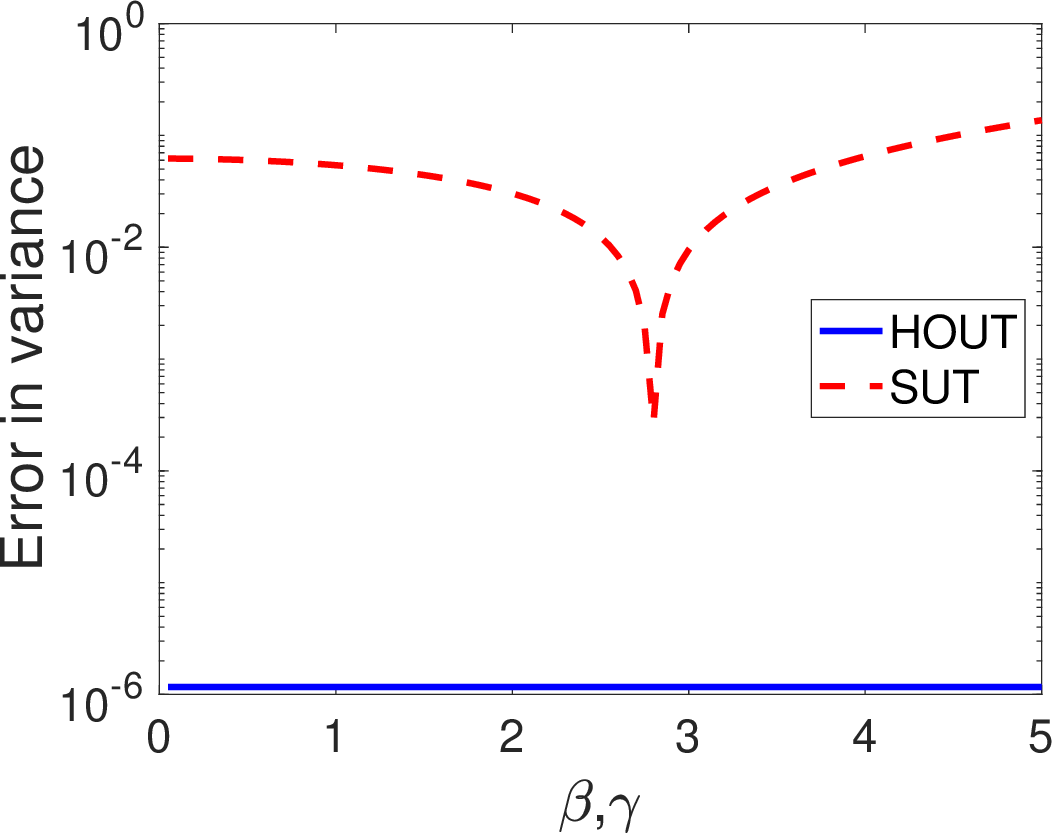}}
   \caption{(a) Comparison between the Higher Order Unscented Transform ensemble (HOUT, red dots) and the Scaled Unscented Transform ensemble (SUT, green dots) on a non-Gaussian distribution. (b,c) Estimating the output mean and covariance for various values of $\beta$ in the SUT and various values of $\gamma$ in the HOUT.}\label{fig2}
  \end{figure}
  
  We first compare the HOUT and SUT on various polynomials applied to a two dimensional input distribution.  In order to generate a non-Gaussian input distribution, we start by generating an ensemble of $10^5$ standard Gaussian random variables, $Z \in \mathbb{R}^2$ and then transforming them by a map $X = AZ + B(Z\odot Z \odot \textup{sign}(Z))$ where $A,B$ are random $2\times 2$ matrices and $\odot$ is componentwise multiplication.  The resulting ensemble is shown in Fig.~\ref{fig2}(a) along with the HOUT (red dots) and SUT (green dots) ensembles.

The SUT has the free parameter $\beta$ but the HOUT requires a certain inequality for $\beta$ and instead the HOUT has $\gamma$ as a free parameter.  In order to explore the effect of these parameters on the SUT and HOUT, we considered a random quadratic polynomial $f:\mathbb{R}^2\to\mathbb{R}$.  In Fig.~\ref{fig2} we show the error of the HOUT and SUT estimates of the mean $\mathbb{E}[f(X)]$ and variance $\mathbb{E}[(f(X)-\mathbb{E}[f(X)])^2]$ as a function of $\beta$ for the SUT and $\gamma$ for the HOUT.  Notice that since $f$ is a quadratic polynomial, the mean is also a quadratic polynomial, whereas the variance is a quartic polynomial.  Since the SUT has degree of exactness two, it is exact on the mean but not on the variance.  The HOUT has degree of exactness four and is exact on both up to the specified tolerance ($10^{-5}$ in these experiments).  Reducing the tolerance below this point led to increased error, most likely due to the conditioning of the HOUT quadrature rule.

\begin{figure}[H]
  \centering
  \subfloat[$f(x) = ax + bcx^2$]
  {\includegraphics[width=0.245\textwidth]{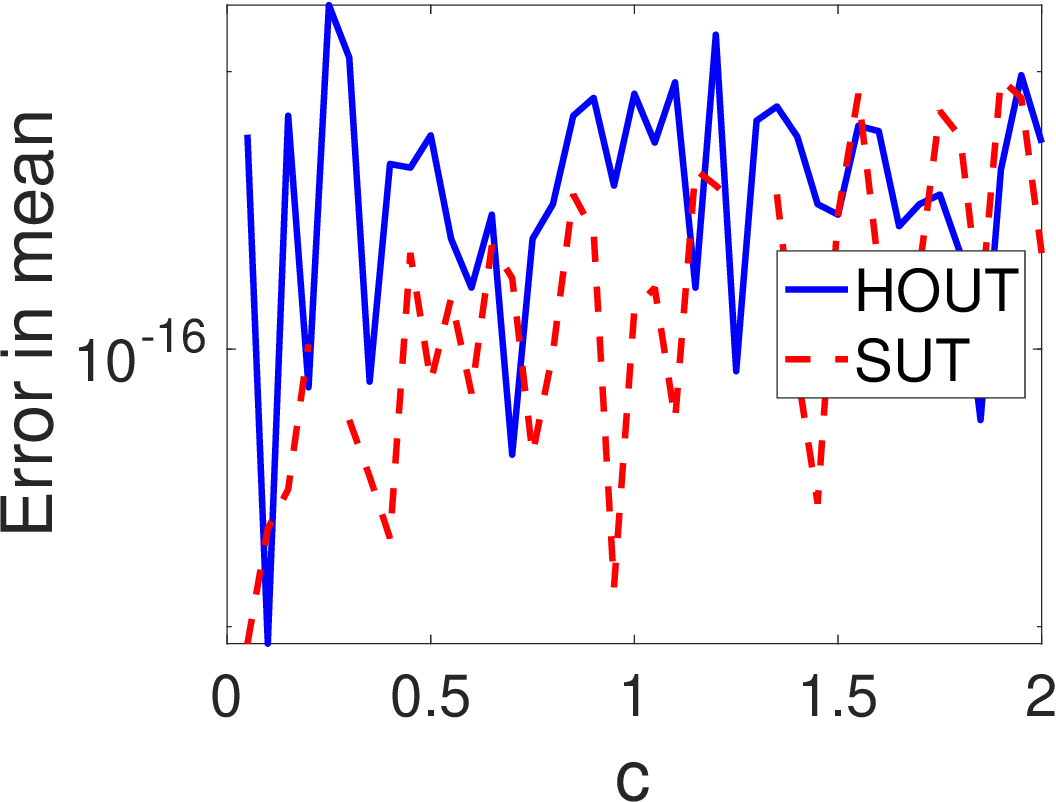}\label{fig:f1}}
  \hfill
  \subfloat[$f(x) = ax + bcx^3$]
  {\includegraphics[width=0.245\textwidth]{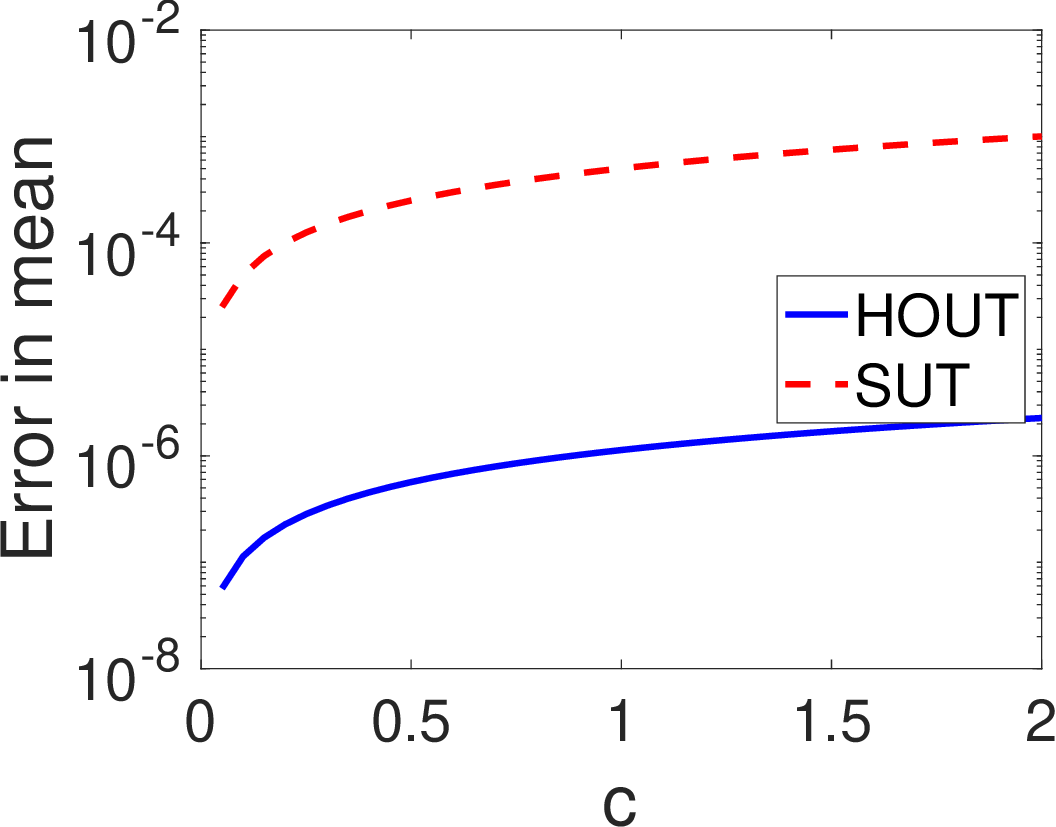}\label{fig:f2}}
  \hfill
  \subfloat[$f(x) = ax + bcx^4$]
  {\includegraphics[width=0.245\textwidth]{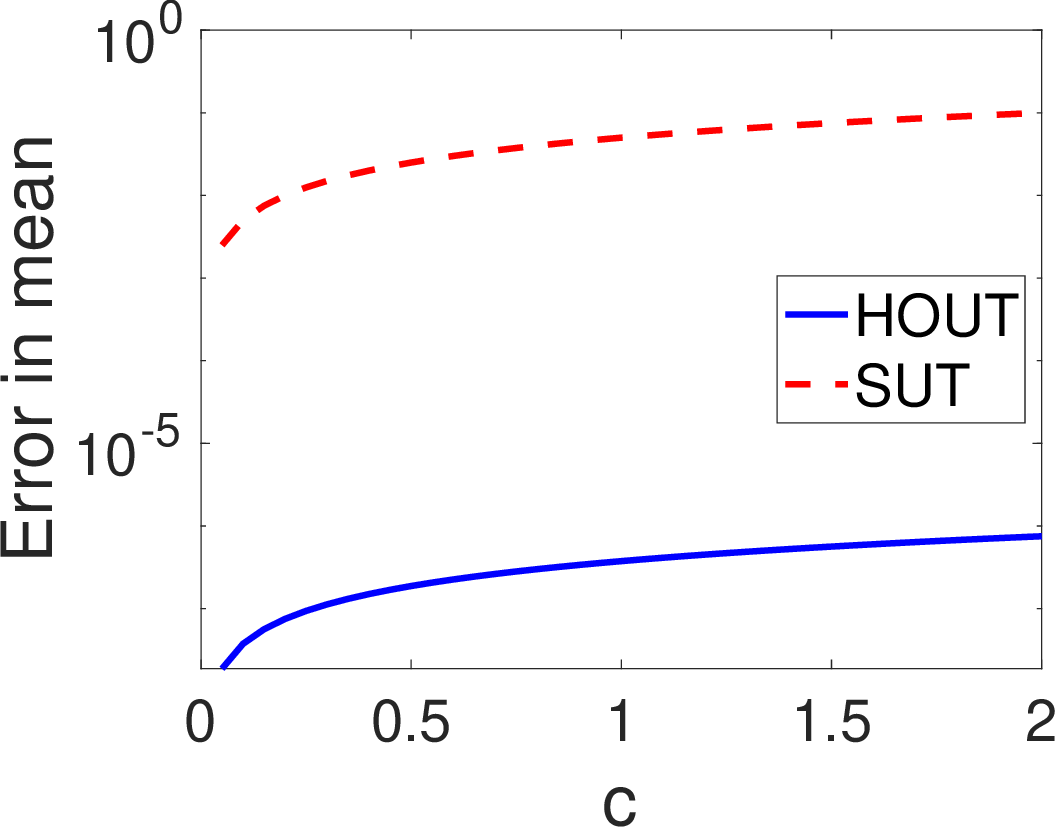}\label{fig:f3}}
  \hfill
  \subfloat[$f(x) = ax + bcx^5$]
  {\includegraphics[width=0.245\textwidth]{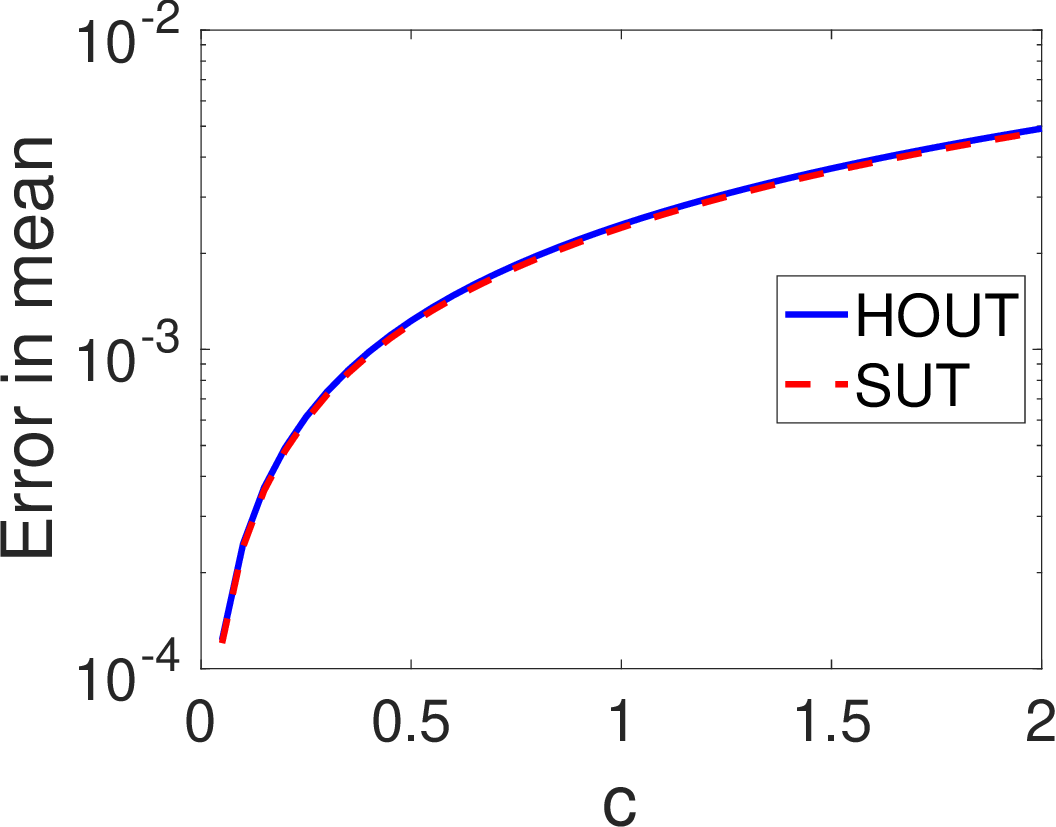}\label{fig:f4}} \\
    \subfloat[$f(x) = ax + bcx^2$]
  {\includegraphics[width=0.245\textwidth]{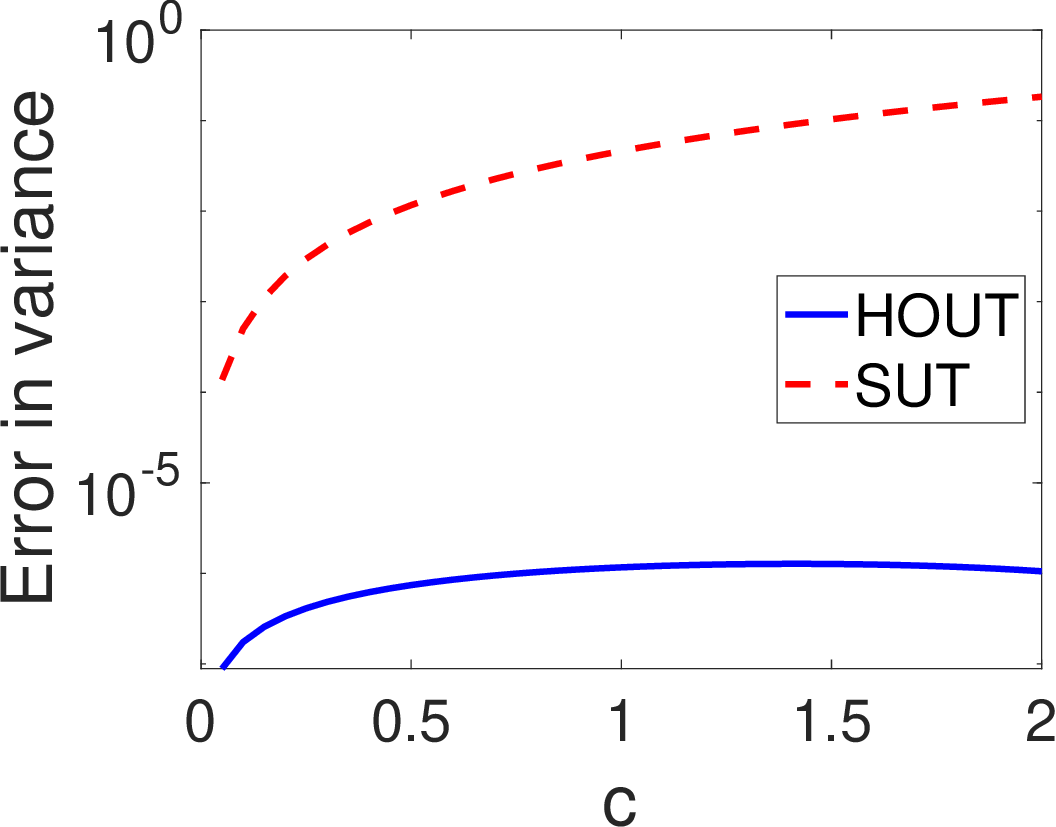}\label{fig:f5}}
  \hfill
  \subfloat[$f(x) = ax + bcx^3$]
  {\includegraphics[width=0.245\textwidth]{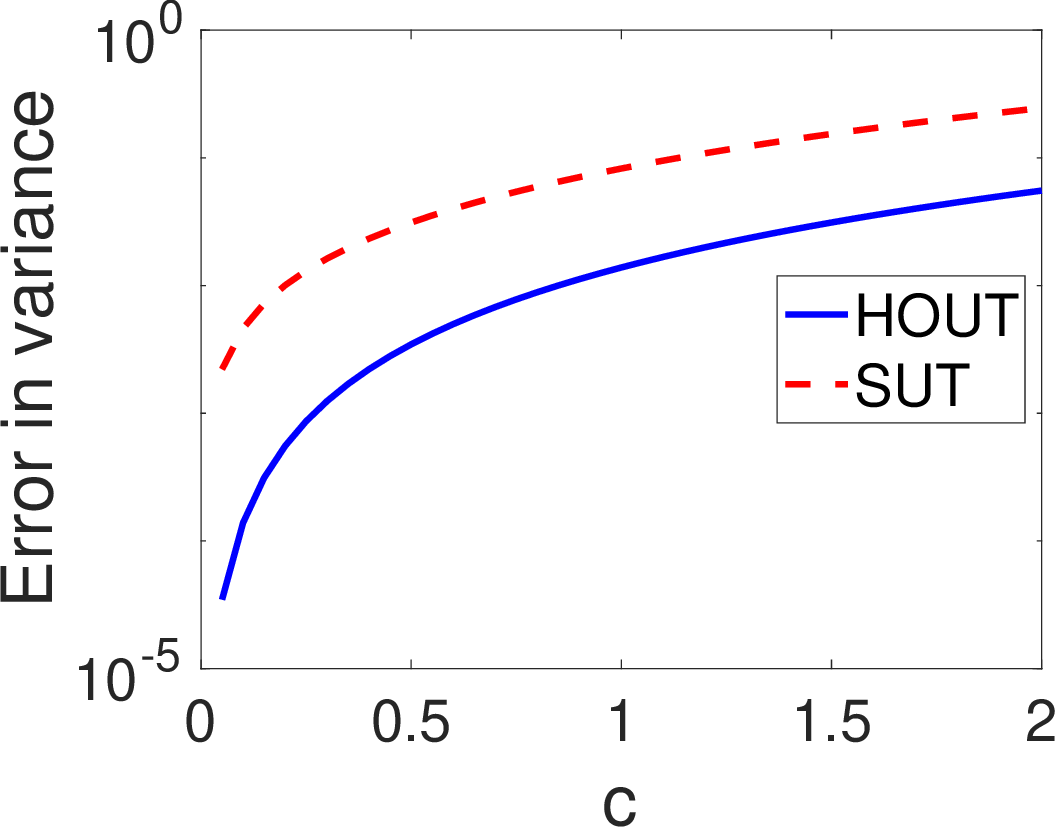}\label{fig:f6}}
  \hfill
  \subfloat[$f(x) = ax + bcx^4$]
  {\includegraphics[width=0.245\textwidth]{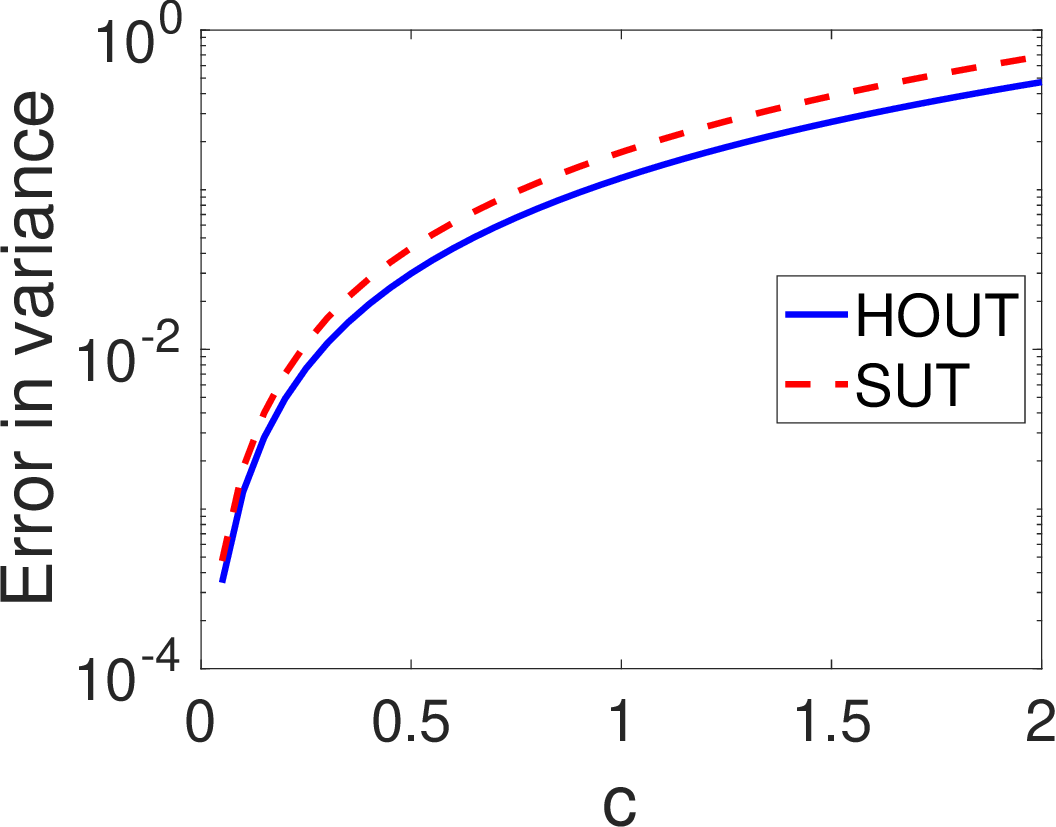}\label{fig:f7}}
  \hfill
  \subfloat[$f(x) = ax + bcx^5$]
  {\includegraphics[width=0.245\textwidth]{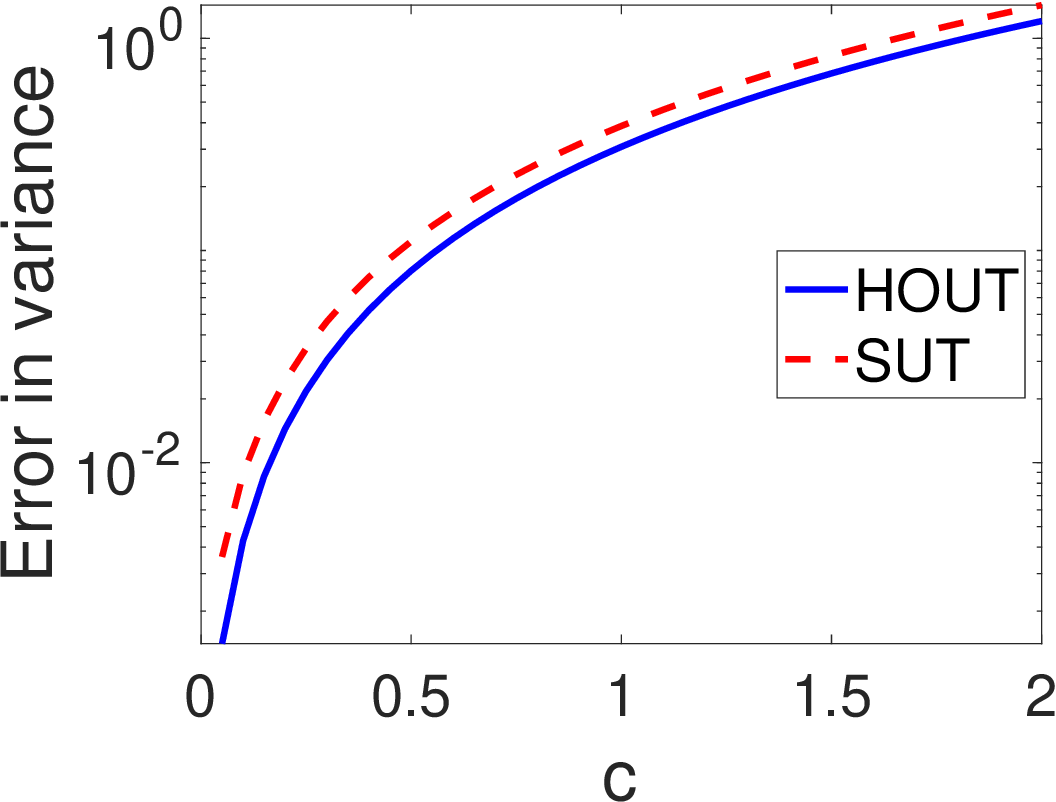}\label{fig:f8}}
  \caption{Comparison between the Higher Order Unscented Transform (HOUT) and the Scaled Unscented Transform (SUT) when estimating the mean $\mathbb{E}[f(X)]$ (top row) and variance $\mathbb{E}[(f(X)-\mathbb{E}[f(X)])^2]$ (bottom row) with different polynomials. Notice that the SUT has degree of exactness two while the HOUT has degree of exactness four.}\label{fig3}
\end{figure}

Using the same two-dimensional distribution, $X$, we passed it through several polynomial functions of the form $f(x) = ax + bcx^n$ for $n=2,3,4,5$ where $a$ and $b$ are made random $1\times 2$ vectors.  To show the influence of the strength of the nonlinearity, we sweep through different values of $c$. In Fig.~\ref{fig3} we compare the HOUT and SUT for estimating the mean and variance of the output of each of these polynomials.  As expected, the HOUT is exact for the means up to $n=4$ and for the variances up to $n=2$ due to having degree of exactness four.  For higher degree polynomials, the HOUT has comparable or better performance.

\begin{figure}[H]
  \centering
  \subfloat[]
  {\includegraphics[width=0.321\textwidth]{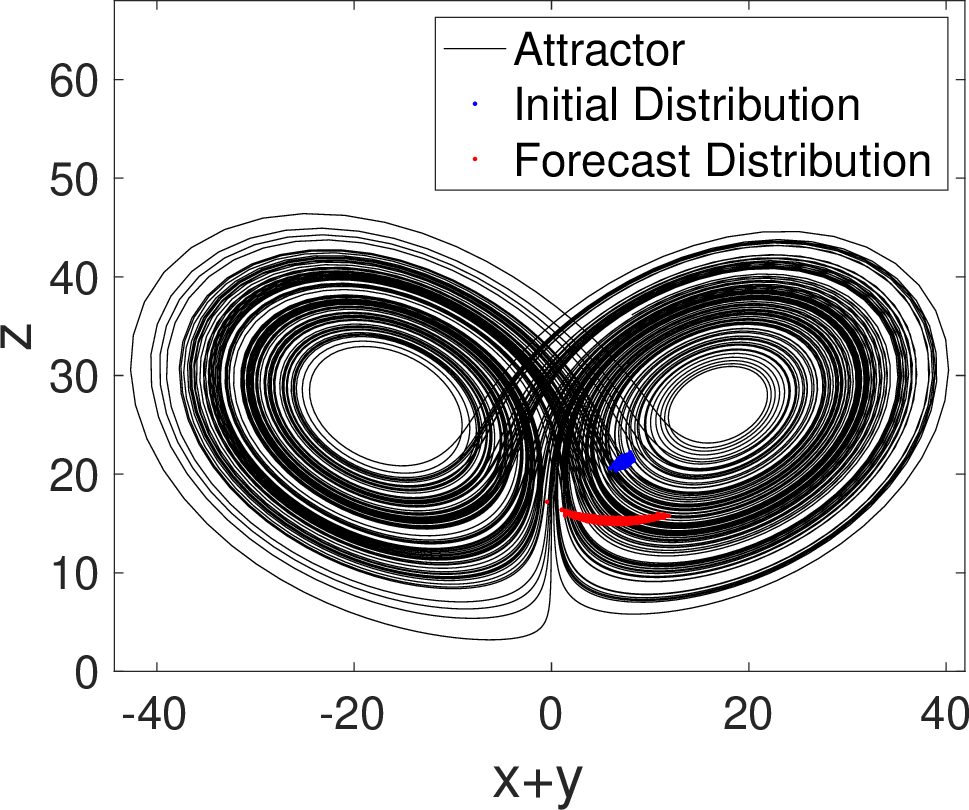}}
  \hfill
  \subfloat[]
  {\includegraphics[width=0.335\textwidth]{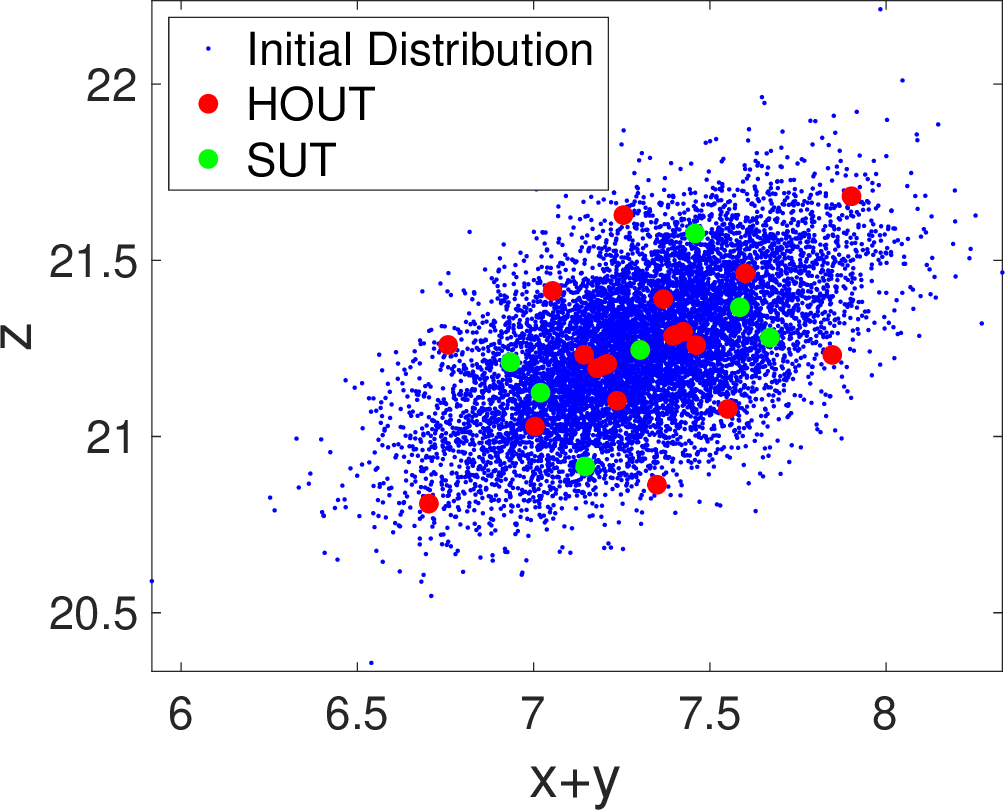}}
  \hfill
  \subfloat[]
  {\includegraphics[width=0.321\textwidth]{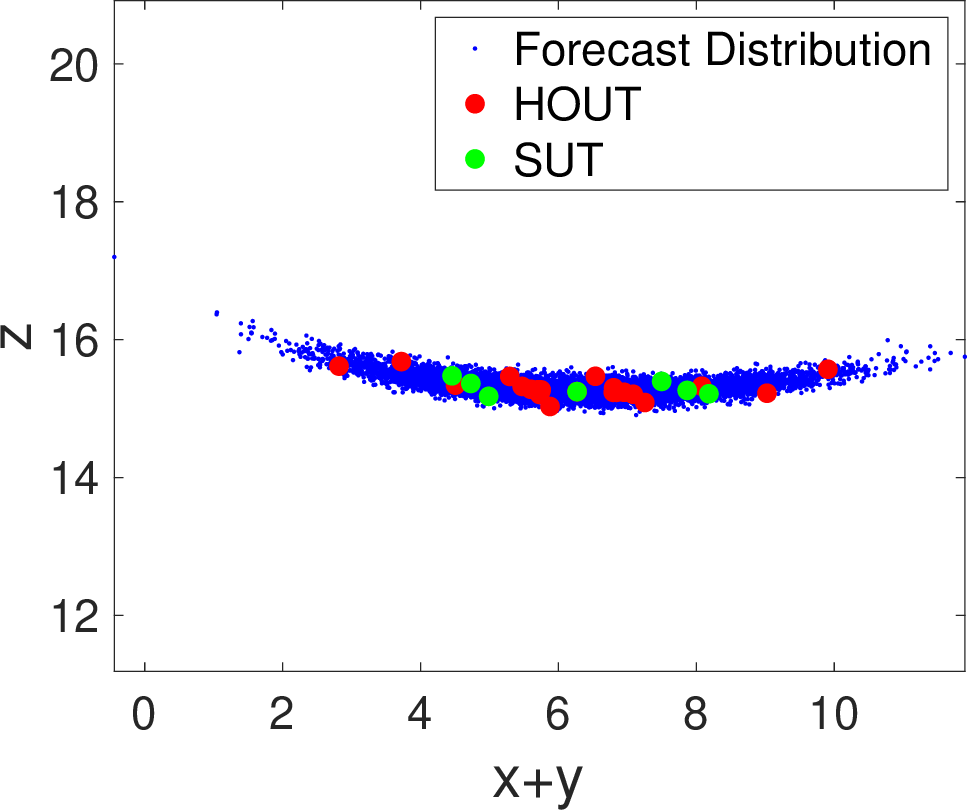}} \\
  \subfloat[]
  {\includegraphics[width=0.245\textwidth]{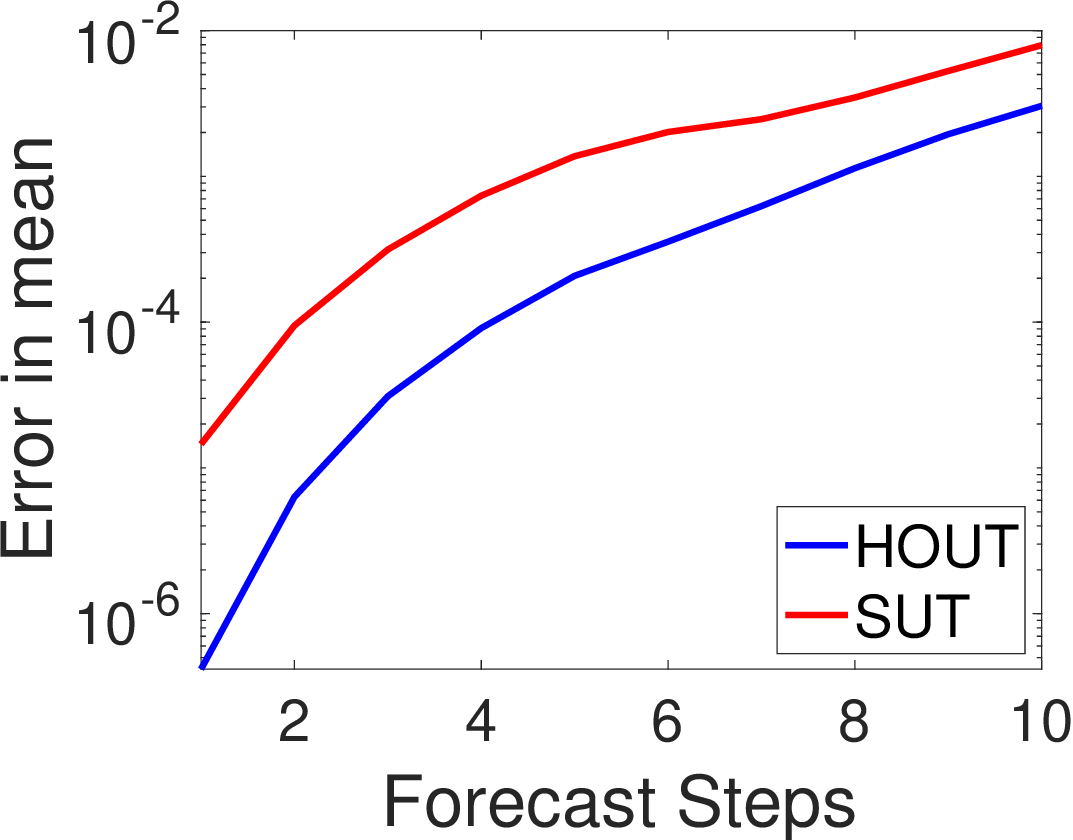}}
  \hfill
  \subfloat[]
  {\includegraphics[width=0.245\textwidth]{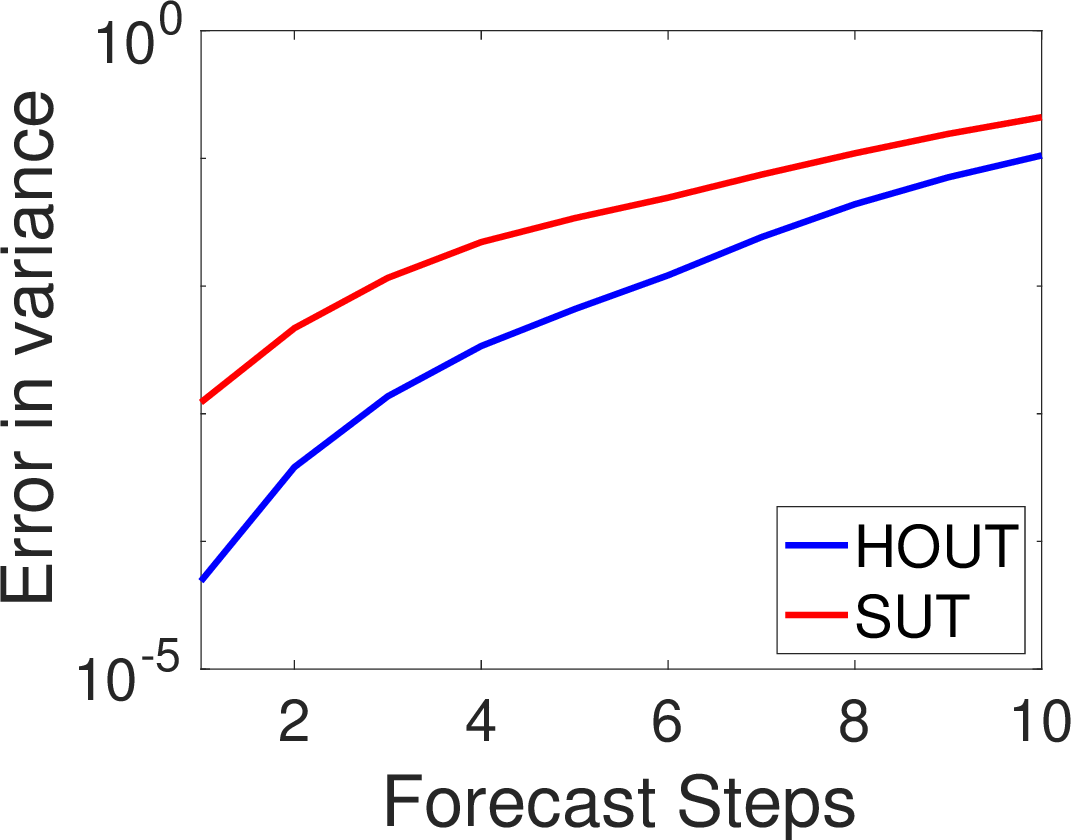}}
  \hfill
  \subfloat[]
  {\includegraphics[width=0.245\textwidth]{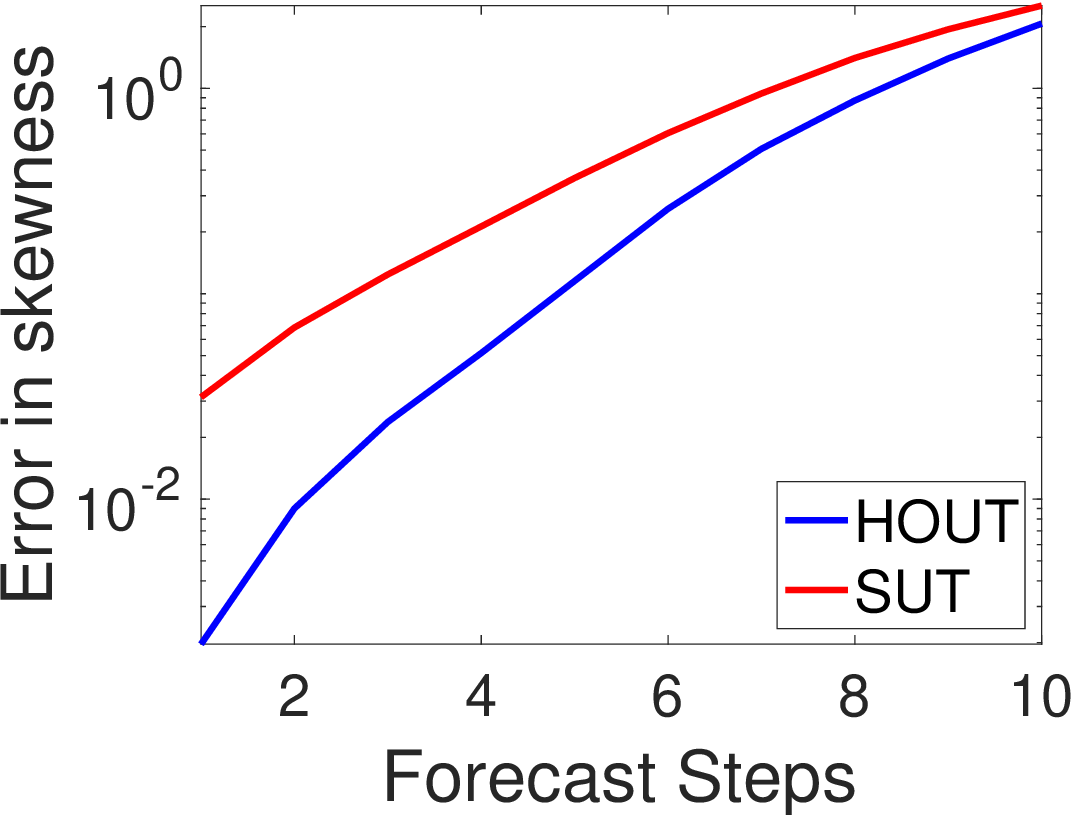}}
  \hfill
  \subfloat[]
  {\includegraphics[width=0.245\textwidth]{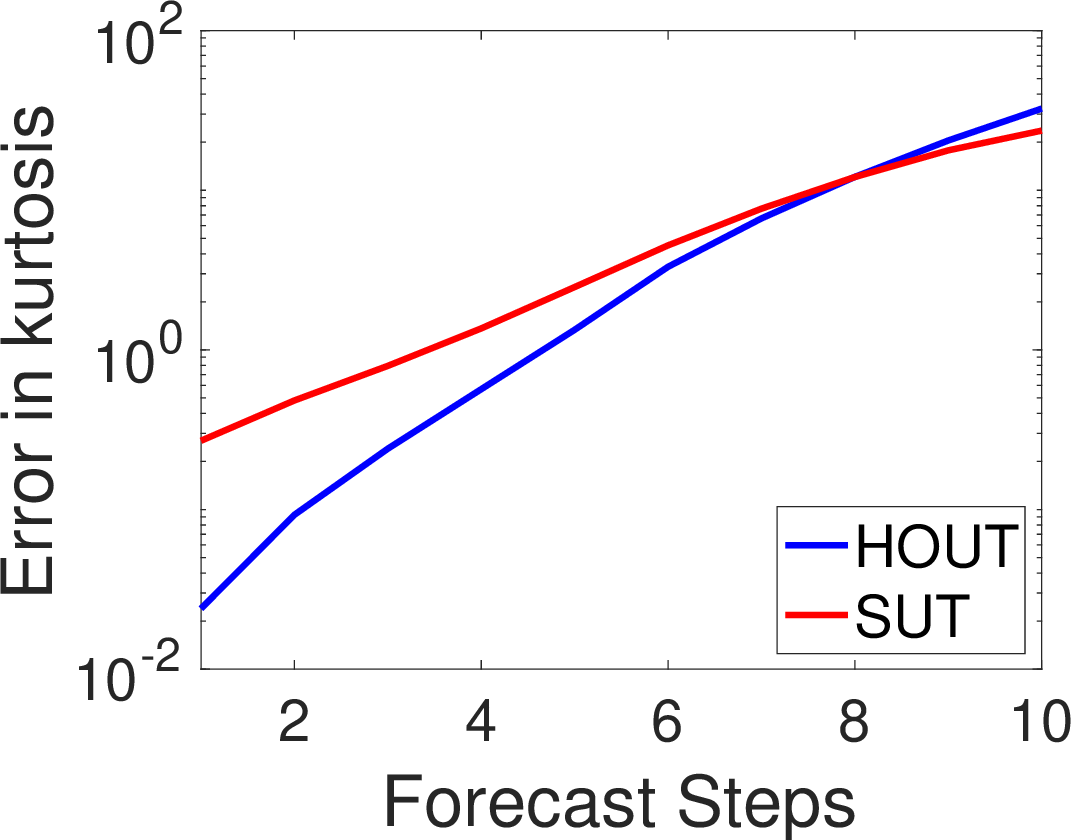}}
  \caption{Comparison between the Higher Order Unscented Transform (HOUT) and the Scaled Unscented Transform (SUT) when estimating the mean $\mathbb{E}[f(X)]$ (top row) and higher moments of the Lorenz-63 model at various forecast horizons.  In (a) we show the Lorenz-63 attractor (black) along with an example initial ensemble (blue) and forecast ensemble (red) used to compute the true statistics.  In (b,c) we show the initial and forecast ensembles (blue) together with the HOUT (red) and SUT (green) ensembles.  Results in (d-g) show the forecast accuracy versus the forecast horizon and are geometrically averaged over 500 different initial conditions on the attractor.}\label{fig4}
\end{figure}

Of course, the HOUT and SUT are intended for use beyond polynomial functions.  In fact, the most common application is for forecasting dynamical systems.  Next, we consider the problem of forecasting the chaotic Lorenz-63 dynamical system \cite{lorenz}.  We integrate the Lorenz-63 system with a Runge-Kutta order four method and a time step $\tau=0.1$.  In order to generate a non-Gaussian initial state, we start by choosing a random point on the attractor and adding a small amount of Gaussian noise.  We then run the ensemble forward $N_1=5$ steps and we consider this the initial state, see Fig.~\ref{fig4}(a) (blue) and Fig.~\ref{fig4}(b) (blue).  We compute the statistics of the initial state using the ensemble shown, and use these statistics to generate the HOUT and SUT as shown in Fig.~\ref{fig4}(b).  All three ensembles are then integrated forward in time $N_2$ additional steps and the true forecast statistics from the large ensemble are compared to the HOUT and SUT estimates.  An example is shown in Fig.~\ref{fig4}(c) with $N_2 = 15$.  

We then repeat this experiment $500$ times with different randomly selected initial states on the attractor and we compute the geometric average of the error between the HOUT estimate and the true statistics at each forecast time, shown in Fig.~\ref{fig4}(d-g)(blue).  Similarly, we compute the geometric average of the error between the SUT estimate and the true statistics (red) at each forecast time, shown in Fig.~\ref{fig4}(d-g)(red).  We note that the HOUT provides improved estimates of the first four moments up to at least $4$ forecast steps, which is $0.4$ model time units.  In particular, the mean forecast is improved by an order of magnitude in this forecast range. 

\section{Conclusions and Future Directions}
\label{sec:conclusions}

The SUT is a highly efficient and successful strategy for uncertainty quantification with many applications.  As computational resources expand, there is a growing demand for larger ensembles that are similarly well designed.  At the same time, complex systems demand better uncertainty quantification such as skewness and kurtosis to capture fat-tails.  The HOUT generalizes the SUT to efficiently leverage additional computation resources to meet the growing UQ demand. There are several promising directions of future research that we expect to result from this work. 

First, there are many applications, such as Kalman filtering and smoothing for nonlinear systems that use the SUT to track the mean and covariance of hidden variables based on noisy observations.  If these filters and smoothers can be generalized to track four moments, they could be integrated with the HOUT to achieve better stability and accuracy along with additional uncertainty quantification.  Moreover, these methods are often difficult to analyze theoretically due to the lack of a natural limit.  This compares to the relative ease of theoretical analysis of particle filters where one may consider the infinite  particle (Monte-Carlo) limit.  The HOUT opens up the possibility that generalized Kalman-based approaches may be analyzed in the limit of infinitely many moments, in a sense this is a kind of spectral solver convergence. While we only explicitly derive the four moment version of the HOUT, the methods used should allow generalization to match an arbitrary number of moments.

A second promising direction for future research concerns deriving error bounds for the SUT and HOUT. Current error bound for the SUT are based on Taylor expansion \cite{julier,julier1995,julier1997}, and a similar analysis could be carried out for the HOUT.  However, this analysis requires decay of the moments and a highly localized input density, moreover it is not the natural method of analyzing quadrature error.  A more natural approach would be based on multivariate polynomial approximation error bounds, which would be analogous to the univariate quadrature error bound analysis.  

Finally, more efficient rank-1 decomposition can immediately improve the efficiency of the HOUT. Similarly, improved/sharp bounds on the relationship between tensor eigenvalues and their entries could improve understanding of the convergence rate as a function of dimension and tensor order.

\appendix


\section{Proof of Lemma \ref{order3and4}} 

\begin{proof}
First note that for 3-tensors, if $\lambda$ is an eigenvalue then $(T\times_1 v)\times_1 v = \lambda v$ so $(T\times_1 (-v))\times_1(-v)= -\lambda(-v)$ so $-\lambda$ is also an eigenvalue. Therefore for 3-tensors, $\lambda_{max}=\lambda_{maxabs}$, and in fact this is true for any odd order tensor.

Next, by the symmetry of the matrix 3-tensor $T$ 
\[\sum_{i,j,k}T_{ijk} v_iv_jv_k = \sum_{i=1}^nT_{iii} v_i^3+3\sum_{i=1}^n\sum_{k\ne i}T_{iik}v_i^2v_k+\sum_{i=1}^n\sum_{j\neq i}\sum_{k\neq i,j} T_{ijk}v_iv_jv_k\]
Let us fix $s=1,...,n$ and let $(v_s)_i=\hbox{sign}(T_{sss})\,\delta_{is}.$ Then $\|v_s\|=1$ so 
\ben*
\l_{max}&=&\max_{\|v\|=1}\sum_{i,j,k}T_{ijk}v_iv_jv_k
\ge \sum_{i,j,k}T_{ijk}(v_s)_i(v_s)_j(v_s)_k\\
&=& (\hbox{sign}(T_{sss}))^3\sum_{i,j,k} T_{ijk}\delta_{is}\delta_{js}\delta_{ks} = \hbox{sign}(T_{sss})\sum_{i,j,k} T_{ijk}\delta_{is}\delta_{js}\delta_{ks}\\
&=& \hbox{sign}(T_{sss})T_{sss} = |T_{sss}|
\eeqn*
Thus, in this case we have $\l_{max}\ge|T_{sss}|$ for all $s=1,...,n$. Next, fix $s,t\in\{1,...,n\}$ and let 
\[(w_{s,t})_i = \hbox{sign}(T_{stt})\frac{\delta_{is}+\delta_{it}}{\sqrt{2}}\]
Then $\|w_{s,t}\|=1$ and so
\ben*
\l_{max} &\ge& \sum_{i,j,k}T_{ijk} (w_{s,t})_i(w_{s,t})_j(w_{s,t})_k = \left(\frac{\hbox{sign}(T_{stt})}{\sqrt{2}}\right)^3(T_{sss}+T_{ttt}+3T_{sst}+3T_{stt})
\eeqn*
and therefore
\begin{eqnarray}\label{appeqn1}
2^{\, 3/2}\l_{max} &\ge& \hbox{sign}(T_{stt})(T_{sss}+T_{ttt}+3T_{sst}+3T_{stt}).
\end{eqnarray}
Now let 
\[(\tilde{w}_{s,t})_i = \hbox{sign}(T_{stt})\frac{\delta_{is}-\delta_{it}}{\sqrt{2}}\] so we have
\ben*
\l_{max} &\ge& \sum_{i,j,k}T_{ijk} (\tilde{w}_{s,t})_i(\tilde{w}_{s,t})_j(\tilde{w}_{s,t})_k = \left(\frac{\hbox{sign}(T_{stt})}{\sqrt{2}}\right)^3(T_{sss}-T_{ttt}-3T_{sst}+3T_{stt})
\eeqn*
and 
\begin{eqnarray}\label{appeqn2}
2^{\, 3/2}\l_{max} &\ge& \hbox{sign}(T_{stt})(T_{sss}-T_{ttt}-3T_{sst}+3T_{stt})
\end{eqnarray}
Adding equations \eqref{appeqn1} and \eqref{appeqn2}, we get 
\ben*
2(2^{\, 3/2}\l_{max}) &\ge& \hbox{sign}(T_{stt})(2T_{sss}+6T_{stt})\\
2^{\, 3/2}\l_{max} &\ge& \hbox{sign}(T_{stt})(T_{sss}+3T_{stt})
\eeqn*
Recall that $\l_{max}\ge|T_{sss}|\ge -\hbox{sign}(T_{stt})T_{sss}$, so
\ben*
2^{\, 3/2}\l_{max} + \l_{max}&\ge& \hbox{sign}(T_{stt})(T_{sss}+3T_{stt}-T_{sss})\\
(2^{\, 3/2}+1)\l_{max} &\ge& 3\,\hbox{sign}(T_{stt})T_{stt}\\
\eeqn*
Therefore \[\l_{max}\ge \frac3{2^{\, 3/2}+1}|T_{stt}|.\]


Lastly, fix $s,t,u\in\{1,...,n\}$ and let 
\[(w_{s,t,u})_i = \hbox{sign}(T_{stu})\frac{\delta_{is}+\delta_{it}+\delta_{iu}}{\sqrt{3}}\]
Then $\|w_{s,t,u}\|=1$ and so
\ben*
\l_{max} \ge \sum_{i,j,k}T_{ijk} (w_{s,t,u})_i(w_{s,t,u})_j(w_{s,t,u})_k 
&=& \left(\frac{\hbox{sign}(T_{stu})}{\sqrt{3}}\right)^3(T_{sss}+T_{ttt}+T_{uuu}+3T_{sst} \\ && +3T_{stt}+3T_{ssu}+3T_{ttu}+3T_{suu}+3T_{tuu}+6T_{stu})
\eeqn*
Note that since $\l_{max}$ is greater than or equal to $-\hbox{sign}(T_{stu})T_{sss}, -\hbox{sign}(T_{stu})T_{ttt},$ and \\ $-\hbox{sign}(T_{stu})T_{uuu}$ we can factor out $\hbox{sign}(T_{stu})$ so that
\ben*
\left(3^{\, 3/2} + 3\right)\l_{max} &\ge& 3\,\hbox{sign}(T_{stu})(T_{sst}+T_{stt}+T_{ssu}+T_{ttu}+T_{suu}+T_{tuu}+2T_{stu})
\eeqn*
Now note that $(2^{\, 3/2}+1)\l_{max}\ge 3|T_{stt}|\ge -3\,\hbox{sign}(T_{stu})T_{stt}$ which implies 
\ben* (2^{\, 3/2}+1)\l_{max}&\ge& -3\,\hbox{sign}(T_{stu})T_{ssu} \\ 
(2^{\, 3/2}+1)\l_{max}&\ge& -3\,\hbox{sign}(T_{stu})T_{ttu} \\ 
(2^{\, 3/2}+1)\l_{max}&\ge& -3\,\hbox{sign}(T_{stu})T_{suu} \\
(2^{\, 3/2}+1)\l_{max}&\ge& -3\,\hbox{sign}(T_{stu})T_{tuu} 
\eeqn*
and together these imply
\ben*
\left(3^{\, 3/2} + 3 + 6(2^{\, 3/2}+1)\right)\l_{max} &\ge& 6\,\hbox{sign}(T_{stu})T_{stu}\\
\left(3 + 4 \sqrt{2} + \sqrt{3}\right)\l_{max} &\ge& 2|T_{stu}|\\
\l_{max} &\ge& \frac{2}{3 + 4 \sqrt{2} + \sqrt{3}} |T_{stu}|
\eeqn*
Comparing the lower bounds found in the above three cases, we see that the conclusion holds if we set $$c=
\frac{2}{3 + 4 \sqrt{2} + \sqrt{3}}.$$ This completes the proof for 3-tensors. Next we follow a similar strategy for 4-tensors.\vspace{1em}

By the symmetry of the matrix 4-tensor $T$, we have 
\ben*
\sum_{i,j,k,\ell}T_{ijk\ell} v_iv_jv_kv_\ell &=& \sum_{i=1}^nT_{iiii} v_i^4+6\sum_{i=1}^n\sum_{j\ne i}T_{iijj}v_i^2v_j^2+4\sum_{i=1}^n\sum_{j\ne i}T_{iiij}v_i^3v_j+12\sum_{i=1}^n\sum_{j\ne i}\sum_{k\ne i,j}T_{iijk}v_i^2v_jv_k\\
&&+24\sum_{i=1}^n\sum_{j\neq i}\sum_{k\neq i,j} \sum_{\ell\neq i,j,k} T_{ijk\ell}v_iv_jv_kv_\ell.
\eeqn*
We wish to show that for some constant $c\in(0,1]$, 
$$\l_{maxabs}\ge c|T_{ijk\ell}|\quad\text{ for all }\quad i,j,k,\ell\in\{1,2,\dots,n\}.$$ We are going to carry out the proof in five steps by looking at the following cases:
\begin{enumerate}
\item 
$i=j=k=\ell$         
\item 
$i=j\neq k=\ell$
\item 
$i=j=k\neq \ell$
\item 
$i=j$ distinct from $k,\ell$ and $k\neq\ell$
\item 
$i,j,k,\ell$ all distinct
\end{enumerate}

\begin{enumerate}
\item Let us fix $s=1,...,n$ and let $(v_s)_i=\delta_{is}.$ Then $\|v_s\|=1$ so 
\ben*
\l_{max}&=&\max_{\|v\|=1}\sum_{i,j,k,\ell}T_{ijk\ell}v_iv_jv_kv_\ell \ge \sum_{i,j,k,\ell}T_{ijk\ell}(v_s)_i(v_s)_j(v_s)_k(v_s)_\ell = T_{ssss}\\
\l_{min}&=&\min_{\|v\|=1}\sum_{i,j,k,\ell}T_{ijk\ell}v_iv_jv_kv_\ell \le \sum_{i,j,k,\ell}T_{ijk\ell}(v_s)_i(v_s)_j(v_s)_k(v_s)_\ell = T_{ssss}.
\eeqn*
Thus $-\lambda_{min}\ge-T_{ssss}$. Therefore, for all $s=1,...,n$, we have 
\ben*\l_{maxabs} &=& \max\{|\l_{max}|,|\l_{min}|\}\ge \l_{max}\ge T_{ssss}\\
\l_{maxabs} &=& \max\{|\l_{max}|,|\l_{min}|\}\ge 
-\l_{min}\ge -T_{ssss}.\eeqn*
Thus \[\l_{maxabs} \ge |T_{ssss}| \quad\text{for each}\quad s=1,...,n.\]
\item Next, fix $s,t\in\{1,...,n\}$ and let 
\[(w_{s,t})_i = \frac{\delta_{is}+\delta_{it}}{\sqrt{2}}.\]
Then $\|w_{s,t}\|=1$ and so
\ben*
\l_{max} &\ge& \sum_{i,j,k,\ell}T_{ijk\ell} (w_{s,t})_i(w_{s,t})_j(w_{s,t})_k(w_{s,t})_\ell\\
&=& \left(\frac1{\sqrt{2}}\right)^4\left(T_{ssss}+T_{tttt}+\binom{4}{2}T_{sstt}+\binom{4}{1}T_{ssst}+\binom{4}{1}T_{sttt}\right).
\eeqn*
Hence
\begin{eqnarray}\label{appeqn3}
4\l_{max} &\ge& T_{ssss}+T_{tttt}+6T_{sstt}+4T_{ssst}+4T_{sttt}.
\end{eqnarray}
Now let 
\[(\tilde{w}_{s,t})_i = \frac{\delta_{is}-\delta_{it}}{\sqrt{2}}.\] Then
\ben*
\l_{max} &\ge& \sum_{i,j,k,\ell}T_{ijk\ell} (\tilde{w}_{s,t})_i(\tilde{w}_{s,t})_j(\tilde{w}_{s,t})_k(\tilde{w}_{s,t})_\ell\\
&=& \left(\frac1{\sqrt{2}}\right)^4\left(T_{ssss}+T_{tttt}+\binom{4}{2}T_{sstt}-\binom{4}{1}T_{ssst}-\binom{4}{1}T_{sttt}\right).
\eeqn*
Thus 
\begin{eqnarray}\label{appeqn4}
4\l_{max} &\ge& T_{ssss}+T_{tttt}+6T_{sstt}-4T_{ssst}-4T_{sttt}
\end{eqnarray}
Adding inequalities \eqref{appeqn3} and \eqref{appeqn4}, and dividing by 2, we obtain
\begin{eqnarray}\label{eq3}
4\l_{max} &\ge& T_{ssss}+T_{tttt}+6T_{sstt}.\end{eqnarray}
Since $-\l_{min}$ is no smaller than both $-T_{ssss}$ and $-T_{tttt}$, subtracting $2\l_{min}$ from \eqref{eq3} gives us \begin{eqnarray}\label{eq4}
6\l_{maxabs}\ge 4\l_{max}-2\l_{min}\ge 6 T_{sstt}
\end{eqnarray}
Therefore 
\begin{eqnarray} \label{lest} \l_{maxabs}\ge T_{sstt}.\end{eqnarray}
Trguing similarly, we can show that 
\begin{eqnarray}\label{eq3.5} 4\l_{min} \le T_{ssss}+T_{tttt}+6T_{sstt}.\end{eqnarray}
Since $-\l_{max}$ is no larger than both $-T_{ssss}$ and $-T_{tttt}$, subtracting $2\l_{max}$ from \eqref{eq3.5} yields 
\begin{eqnarray} \label{minmaxsstt} 4\l_{min}-2\l_{max}\le 6 T_{sstt}.\end{eqnarray}
Thus 
$$6\l_{maxabs}\ge 2\l_{max}-4\l_{min}\ge -6 T_{sstt}.$$

\noindent Hence $\l_{maxabs}\ge -T_{sstt}$. Using this and \eqref{lest}, we obtain
\[\l_{maxabs}\ge |T_{sstt}|.\]
\item Next, fix distinct $s,t\in\{1,...,n\}$ and let 
\[(w_{s,t})_i = \frac{\delta_{is}-2\delta_{it}}{\sqrt{3}}.\]
Then $\|w_{s,t}\|=1$ and so
\ben*
\l_{max} &\ge& \sum_{i,j,k,\ell}T_{ijk\ell} (w_{s,t})_i(w_{s,t})_j(w_{s,t})_k(w_{s,t})_\ell\\
&=& \left(\frac1{\sqrt{3}}\right)^4\left(T_{ssss}+16T_{tttt}+24T_{sstt}-8T_{ssst}-32T_{sttt}\right).
\eeqn*
Hence
\begin{eqnarray}\label{eq5}
9\l_{max} &\ge& T_{ssss}+16T_{tttt}+24T_{sstt}-8T_{ssst}-32T_{sttt}.
\end{eqnarray}
Adding inequalities \eqref{appeqn3} multiplied by 8 and \eqref{eq5}, we get 
\begin{eqnarray}\label{41lmax}
41\l_{max} &\ge& 9T_{ssss}+24T_{tttt}+72T_{sstt}+24T_{ssst}.
\end{eqnarray}
Since $-\l_{min}$ is no smaller than both $-T_{ssss}$ and $-T_{tttt}$, 
\begin{eqnarray}\label{33lmin} -33\l_{min}\ge -9T_{ssss}-24T_{tttt}.\end{eqnarray}
 Moreover, by \eqref{minmaxsstt}, 
\begin{eqnarray}\label{max2min} \l_{max}-2\l_{min}\ge-3T_{sstt}.\end{eqnarray} Thus, by \eqref{max2min} and using \eqref{41lmax} and \eqref{33lmin}, we have
\ben*
65\l_{max} -81\l_{min}&=&41\l_{max} -33\l_{min}+24(\l_{max}-2\l_{min})\\&\ge& 9T_{ssss}+24T_{tttt}+72T_{sstt}+24T_{ssst}\\&&\hskip10pt -\,9T_{ssss}-24T_{tttt}-72T_{sstt}\\
&=& 24T_{ssst}.
\eeqn*
Thus 
\begin{eqnarray}\label{146} 146\l_{maxabs}\ge 65\l_{max} -81\l_{min}\ge 24T_{ssst}.\end{eqnarray}
Similarly, we can show that
\[65\l_{min} -81\l_{max}\le 24T_{ssst}.\]
Thus
\begin{eqnarray}\label{plus146} 146\l_{maxabs}\ge 81\l_{max} -65\l_{min}\ge -24T_{ssst}.\end{eqnarray}\vspace{0.5em}

\noindent Therefore, by \eqref{146} and \eqref{plus146}, we obtain
\[\l_{maxabs}\ge \frac{12}{73}|T_{ssst}|.\]

\item Next, fix distinct $s,t,u\in\{1,...,n\}$ and let \[(w_{s,t,u})_i = \frac{\delta_{is}+\delta_{it}+\delta_{iu}}{\sqrt{3}}.\] Then $\|w_{s,t,u}\|=1$ and so \ben*\l_{max} &\ge& \sum_{i,j,k,\ell}T_{ijk\ell} (w_{s,t,u})_i(w_{s,t,u})_j(w_{s,t,u})_k(w_{s,t,u})_\ell\\&=& \left(\frac1{\sqrt{3}}\right)^4(T_{ssss}+T_{tttt}+T_{uuuu}+12T_{sstu}+12T_{sttu}+12T_{stuu}+6T_{sstt} \\ &&\hspace{10pt}+6T_{ssuu}+6T_{ttuu}+4T_{ssst}+4T_{sssu}+4T_{sttt}+4T_{suuu}+4T_{tuuu}+4T_{uttt}).\eeqn*
Thus,
\begin{eqnarray}\label{eq7} 9\l_{max}&\ge& T_{ssss}+T_{tttt}+T_{uuuu}+12T_{sstu}+12T_{sttu}+12T_{stuu}\\
&&\hskip10pt +6T_{sstt}+6T_{ssuu}+6T_{ttuu}+4T_{ssst}+4T_{sssu}+4T_{sttt}\nonumber\\
&&\hskip10pt+4T_{suuu}+4T_{tuuu}+4T_{uttt}.\nonumber\end{eqnarray}
Now let 
\[(\tilde{w}_{s,t,u})_i = \frac{\delta_{is}-\delta_{it}-\delta_{iu}}{\sqrt{3}}.\] Then
\ben*
\l_{max} &\ge& \sum_{i,j,k,\ell}T_{ijk\ell} (\tilde{w}_{s,t,u})_i(\tilde{w}_{s,t,u})_j(\tilde{w}_{s,t,u})_k(\tilde{w}_{s,t,u})_\ell\\&=& \left(\frac1{\sqrt{3}}\right)^4(T_{ssss}+T_{tttt}+T_{uuuu}+12T_{sstu}-12T_{sttu}-12T_{stuu}+6T_{sstt}\\ &&\hspace{10pt}+6T_{ssuu}+6T_{ttuu}-4T_{ssst}-4T_{sssu}-4T_{sttt}-4T_{suuu}+4T_{tuuu}+4T_{uttt}).
\eeqn*
Thus,
\begin{eqnarray}\label{eq8} 9\l_{max}&\ge& T_{ssss}+T_{tttt}+T_{uuuu}+12T_{sstu}-12T_{sttu}-12T_{stuu}\\
&& \hskip10pt+6T_{sstt}+6T_{ssuu}+6T_{ttuu}-4T_{ssst}-4T_{sssu}-4T_{sttt}\nonumber\\
&& \hskip10pt-4T_{suuu}+4T_{tuuu}+4T_{uttt}.\nonumber\end{eqnarray}
Adding inequalities \eqref{eq7} and \eqref{eq8}, we get 
\begin{eqnarray}\label{eq9} 18\l_{max}&\ge& 2(T_{ssss}+T_{tttt}+T_{uuuu})+24T_{sstu}\\&& \hskip10pt+12(T_{sstt}+T_{ssuu}+T_{ttuu})+8(T_{tuuu}+T_{uttt}).\nonumber
\end{eqnarray}

Recall from \eqref{plus146} that
\begin{eqnarray} \label{81max} 81\l_{max}-65\l_{min}\ge -24 T_{ssst}.\end{eqnarray} Tpplying twice this estimate to the indices $t$ and $u$ aftering dividing both sides by 3, we obtain
\begin{eqnarray}\label{54} 54\l_{max}-\frac{130}3\l_{min}\ge -8(T_{tuuu}+T_{uttt}).\end{eqnarray}
Since $-\l_{min}$ is no smaller than $-T_{ssss}$, $-T_{tttt}$, and $-T_{uuuu}$, using \eqref{max2min}
 and \eqref{54}, 
inequality \eqref{eq9} yields
\ben*
18\l_{max}&-&6\l_{min}+12\l_{max}-24\l_{min}+54\l_{max}-\frac{130}3\l_{min} \\ &\ge& 2(T_{ssss}+T_{tttt}+T_{uuuu})+24T_{sstu}+12(T_{sstt}+T_{ssuu}+T_{ttuu})\\&&+8(T_{tuuu}+T_{uttt})-2(T_{ssss}+T_{tttt}+T_{uuuu})\\&&-12(T_{sstt}+T_{ssuu}+T_{ttuu})-8(T_{tuuu}+T_{uttt}).\eeqn*
Hence
\ben* 84\l_{max}-\frac{220}3\l_{min}&\ge& 24T_{sstu}.
\eeqn*
Multiplying by 3, we obtain \[472\l_{maxabs}\ge252\l_{max}-220\l_{min}\ge 72T_{sstu}
\]
Following the same argument, we can show that
\[252\l_{min}-220\l_{max}\le 72T_{sstu}.\]
Thus
\begin{eqnarray}\label{220max}472\l_{maxabs}\ge 220\l_{max}-252\l_{min}\ge -72T_{sstu}.\end{eqnarray}Therefore
\[\l_{maxabs}\ge \frac{9}{59}|T_{sstu}|.\]\vspace{1em}

\item Next, fix distinct $s,t,u,v\in\{1,...,n\}$ and let \[(w_{s,t,u,v})_i = \frac{\delta_{is}+\delta_{it}+\delta_{iu}+\delta_{iv}}{2}.\] Then $\|w_{s,t,u,v}\|=1$ and so \ben*\l_{max} &\ge& \sum_{i,j,k,\ell}T_{ijk\ell} (w_{s,t,u,v})_i(w_{s,t,u,v})_j(w_{s,t,u,v})_k(w_{s,t,u,v})_\ell\\&=& \left(\frac12\right)^4(T_{ssss}+T_{tttt}+T_{uuuu}+T_{vvvv}\\
&&\hskip12pt 
+\,6(T_{sstt}+T_{ssuu}+T_{ssvv}+T_{ttuu}+T_{ttvv}+T_{uuvv})\\
&&\hskip12pt +\,12(T_{sstu}+T_{sstv}+T_{ssuv}+T_{ttsu}+T_{ttsv}+T_{ttuv}\\
&&\hskip12pt +\,T_{uust}+T_{uusv}+T_{uutv}+T_{vvst}+T_{vvsu}+T_{vvtu})\\
&&\hskip12pt
+\,4(T_{ssst}+T_{sssu}+T_{sssv}+T_{ttts}+T_{tttu}+T_{tttv}\\
&&\hskip12pt +\,T_{uuus}+T_{uuut}+T_{uuuv}+T_{vvvs}+T_{vvvt}+T_{vvvu})\\
&&\hskip12pt +\,24T_{stuv}.\eeqn*
Thus,
\ben* 16\l_{max}&\ge& T_{ssss}+T_{tttt}+T_{uuuu}+T_{vvvv}\\
&&\hskip12pt +\,6(T_{sstt}+T_{ssuu}+T_{ssvv}+T_{ttuu}+T_{ttvv}+T_{uuvv})\\
&&\hskip12pt +\,12(T_{sstu}+T_{sstv}+T_{ssuv}+T_{ttsu}+T_{ttsv}+T_{ttuv}\\&&\hskip12pt +\,T_{uust}+T_{uusv}+T_{uutv}+T_{vvst}+T_{vvsu}+T_{vvtu})\\
&&\hskip12pt
+\,4(T_{ssst}+T_{sssu}+T_{sssv}+T_{ttts}+T_{tttu}+T_{tttv}\\
&&\hskip12pt +\,T_{uuus}+T_{uuut}+T_{uuuv}+T_{vvvs}+T_{vvvt}+T_{vvvu})\nonumber\\
&&\hskip12pt +\,24T_{stuv}.\nonumber\eeqn*
Since $-\l_{min}$ is no smaller than the quantities $-T_{ssss}$, $-T_{tttt}$, $-T_{uuuu}$ and $-T_{vvvv}$, by \eqref{max2min} applied to the pairs of indices $\{s,t\}$, $\{s,u\}$, $\{s,v\}$, $\{t,u\}$, $\{t,v\}$, and $\{u,v\}$, 
 and, in addition, using \eqref{81max} and \eqref{220max},
 we have
\ben*
16\l_{max}-4\l_{min}+12\l_{max}-24\l_{min}+440\l_{max}-504\l_{min}+162\l_{max}-130\l_{min}&\ge& 24T_{stuv}\eeqn*
which reduces to 
\ben* 
630\l_{max}-662\l_{min}&\ge& 24T_{stuv}.
\eeqn*
Hence \[1292\l_{maxabs}\ge630\l_{max}-662\l_{min}\ge 24T_{stuv}.
\]
Similarly, we can show that
\[630\l_{min}-662\l_{max}\le 24T_{stuv}.\]
Thus
\[1292\l_{maxabs}\ge 662\l_{max}-630\l_{min}\ge -24T_{stuv}.\]Therefore, simplifying, we obtain
\[\l_{maxabs}\ge \frac{6}{323}|T_{stuv}|.\]
Comparing the lower bounds found in the above three cases, we see that the conclusion holds if we set \[c=
\frac{6}{323}.\]

\end{enumerate}

\end{proof}

\section{Proof of Theorem \ref{fourmoments}}

\begin{proof} Suppose $
\eta = \frac1{2}\beta^{-2}$, $\psi = \frac1{2}\rho^{-4}$, $\nu =\frac1{2}\gamma^{-3}$ and $\delta^2 = \rho^2$.
We first we wish to show that the first moment equation matches our mean. We begin by splitting the sum 
\[\sum_{i=-2}^{N} w_i\sigma_i = \sum_{i=-2}^{0} w_i\sigma_i+ \sum_{i=1}^{2d} w_i\sigma_i+\sum_{i=2d+1}^{2d+2J} w_i\sigma_i+ \sum_{i=2d+2J+1}^{N} w_i\sigma_i
\]
Using the expressions defining the 4 moment $\sigma$-points $\sigma_i$ and the corresponding weights $w_i$, we have 
\ben*
\sum_{i=-2}^{N} w_i\sigma_i &=& (1-d\beta^{-2}-\hat L\delta^{-4})\mu + \frac1{2\alpha}(\mu+\alpha\hat{\mu})-\frac1{2\alpha}(\mu-\alpha\hat{\mu}) \\
&+&\sum_{i=1}^d \frac1{2\beta^2}\big(\mu+\beta\sqrt{\hat{C}}_i\big)+\sum_{j=d+1}^{2d}\frac1{2\beta^2}\big(\mu-\beta\sqrt{\hat{C}}_{i-d}\big)\\
&&+\sum_{i=2d+1}^{2d+J}\frac1{2\gamma^3}\big(\mu+\gamma\tilde{v}_{i-2d}\big)+\sum_{j=2d+J+1}^{2d+2J}\frac{-1}{2\gamma^3}\big(\mu-\gamma\tilde{v}_{i-2d-J}\big)\\
&&+\sum_{i=2d+2J+1}^{2d+2J+L}\frac1{2\delta^4} s_{i-2d-2J}\big(\mu+\delta\tilde{u}_{i-2d-2J}\big) \\
&&+\sum_{j=2d+2J+L+1}^{N}\frac1{2\delta^4} s_{i-2d-2J-L}\big(\mu-\delta\tilde{u}_{i-2d-2J-L})\\
\eeqn*
and regrouping like terms, we obtain
\ben*
\sum_{i=-2}^{N} w_i\sigma_i &=& (1-d\beta^{-2}-\hat L\delta^{-4})\mu + \hat{\mu}+\sum_{i=1}^d \frac1{2\beta^2}\big(2\mu+\beta\sqrt{\hat{C}}_i-\beta\sqrt{\hat{C}}_i\big) \\&& 
+\sum_{i=1}^J\left(\frac1{2\gamma^3}\big(\mu+\gamma\tilde{v}_i\big)-\frac1{2\gamma^3}\big(\mu-\gamma\tilde{v}_i\big)\right)+\sum_{i=1}^L\frac1{2\delta^4} s_i\big(2\mu+\delta\tilde{u}_i-\delta\tilde{u}_i\big)\\
&=& (1-d\beta^{-2}-\hat L\delta^{-4})\mu + \hat{\mu}+\sum_{i=1}^d\frac{\mu}{\beta^2}+\sum_{i=1}^J \frac{\tilde{v}_i}{\gamma^2} +\sum_{i=1}^L \frac{s_i\mu}{\delta^4} \\
&=& (1-d\beta^{-2}-\hat L\delta^{-4})\mu+ \hat{\mu}+d\beta^{-2}\mu+\gamma^{-2}\sum_{i=1}^J\tilde{v}_i +\delta^{-4}\mu\sum_{i=1}^Ls_i\\
&=& \mu + \hat{\mu}+\gamma^{-2}\tilde\mu \\
&=& \mu
\eeqn*
using the definition $\hat{\mu} = - \gamma^{-2}\tilde{\mu}$ for the last equality. \vspace{1em}

To look at the other moment equations, let's first observe that  for $n=2,3,4$,
\ben*
\sum_{i=-2}^{N} w_i(\sigma_i-\mu)^{\otimes n} &=& \sum_{i=-1}^{0} w_i(\sigma_i-\mu)^{\otimes n}+ \sum_{i=1}^{2d} w_i(\sigma_i-\mu)^{\otimes n}\\&&+\sum_{i=2d+1}^{2d+2J} w_i(\sigma_i-\mu)^{\otimes n}+ \sum_{i=2d+2J+1}^{N} w_i(\sigma_i-\mu)^{\otimes n}
\eeqn*
By the definition of $\sigma$-points and corresponding weights, 
\ben*
\sum_{i=-2}^{N} w_i(\sigma_i-\mu)^{\otimes n} 
 &=& \frac{\alpha^{n-1}}{2}\left(\hat{\mu}^{\otimes n}-(-\hat{\mu})^{\otimes n}\right)+\frac{\beta^{n-2}}{2}\left(\sum_{i=1}^{d} \big(\sqrt{\hat{C}}_i\big)^{\otimes n}+\sum_{j=d+1}^{2d} \big(-\sqrt{\hat{C}}_{i-d}\big)^{\otimes n}\right)\\
 &&+\frac{\gamma^{n-3}}{2}\left(\sum_{i=2d+1}^{2d+J}\big(\tilde{v}_{i-2d}\big)^{\otimes n}-\sum_{j=2d+J+1}^{2d+2J}\big(-\tilde{v}_{i-2d-J}\big)^{\otimes n} \right) \\
 &&+\frac{\delta^{n-4}}{2}\sum_{i=2d+2J+1}^{2d+2J+L} s_{i-2d-2J}\big(\tilde{u}_{i-2d-2J}\big)^{\otimes n}\\
 &&+\frac{\delta^{n-4}}{2}\sum_{j=2d+2J+L+1}^{N}  s_{i-2d-2J-L}\big(-\tilde{u}_{i-2d-2J-L}\big)^{\otimes n}.
\eeqn*
where we used the property $(av)^{\otimes n}=a^nv^{\otimes n}$ 
where $a$ is any real number and $v$ is a vector.  When $n$ is even, we have 
\begin{equation} \label{eq:even}
\sum_{i=-2}^{N} w_i(\sigma_i-\mu)^{\otimes n} 
 = \beta^{n-2}\sum_{i=1}^{d} \big(\sqrt{\hat{C}}_i\big)^{\otimes n}+\delta^{n-4}\sum_{i=1}^{L} s_{i}\big(\tilde{u}_{i}\big)^{\otimes n}
\end{equation}
and when $n$ is odd, we obtain
\begin{equation} \label{eq:odd}
\sum_{i=-2}^{N} w_i(\sigma_i-\mu)^{\otimes n} 
 = \alpha^{n-1}\hat{\mu}^{\otimes n}+\gamma^{n-3}\sum_{i=1}^{J}\big(\tilde{v}_{i}\big)^{\otimes n}.
\end{equation}

Now we wish to show that the second moment equation matches our covariance. By \eqref{eq:even}, setting $n=2$ we have
\[\sum_{i=-2}^{N} w_i(\sigma_i-\mu)^{\otimes 2} 
  = \sum_{i=1}^{d} \sqrt{\hat{C}}_i^{\otimes 2}+\delta^{-2}\sum_{i=1}^{L}s_i\tilde{u}_i^{\otimes 2} = \hat{C}+\delta^{-2}\tilde{C}\]
  and applying the definition of $\hat{C} = C - \delta^{-2}\tilde C$ we have
  \ben*
  \sum_{i=-2}^{N} w_i(\sigma_i-\mu)^{\otimes 2}  = C,
\eeqn*
as desired.  Next, observe that by \eqref{eq:odd} and the definition of $S$,
\ben*
\sum_{j=0}^{N} w_i(\sigma_i-\mu)^{\otimes 3}
&=& \alpha^2\hat{\mu}^{\otimes 3}+\sum_{i=1}^{J}\tilde{v}_i^{\otimes 3}\\
\eeqn*
and since we assume that $\left|\left|\sum_{i=1}^{J}\tilde{v}_i^{\otimes 3} - S \right|\right|_F \leq \tau/2$, by the triangle inequality we have,
\[ \left|\left|\sum_{j=0}^{N} w_i(\sigma_i-\mu)^{\otimes 3} - S \right|\right|_F \leq \tau/2 + \alpha^2 ||\hat \mu^{\otimes 3}||_F \] 
as desired.  Lastly, we wish to show that the fourth moment equation matches our kurtosis. By \eqref{eq:even} and the definition of $K$, we have
\ben*
\sum_{i=-2}^{N} w_i(\sigma_i-\mu)^{\otimes 4}
  &=& \beta^2\sum_{i=1}^{d} \sqrt{\hat{C}}_i^{\otimes 4}+ \sum_{i=1}^{L}s_i\tilde{u}_i^{\otimes 4} \\
  &=& \beta^2\bar C+ \sum_{i=1}^{L}s_i\tilde{u}_i^{\otimes 4}
\eeqn*
and since we assume that $\left|\left|\sum_{i=1}^{L}s_i\tilde{u}_i^{\otimes 4} - K \right|\right|_F \leq \tau/2$, by the triangle inequality we have,
\[ \left|\left|\sum_{j=0}^{N} w_i(\sigma_i-\mu)^{\otimes 4} - K \right|\right|_F \leq \tau/2 + \beta^2 ||\bar C||_F \] 
which completes the proof.
\end{proof}

\bibliographystyle{siamplain}
\bibliography{references}

\begin{thebibliography}{10}

\bibitem{anderson2007adaptive}
{\sc J.~L. Anderson}, {\em An adaptive covariance inflation error correction
  algorithm for ensemble filters}, Tellus A: Dynamic Meteorology and
  Oceanography, 59 (2007), pp.~210--224.

\bibitem{cubaturefilter}
{\sc I.~Arasaratnam and S.~Haykin}, {\em Cubature kalman filters}, IEEE
  Transactions on automatic control, 54 (2009), pp.~1254--1269.

\bibitem{berry2013adaptive}
{\sc T.~Berry and T.~Sauer}, {\em Adaptive ensemble kalman filtering of
  non-linear systems}, Tellus A: Dynamic Meteorology and Oceanography, 65
  (2013), p.~20331.

\bibitem{HOPMoriginal}
{\sc L.~De~Lathauwer, P.~Comon, B.~De~Moor, and J.~Vandewalle}, {\em
  Higher-order power method}, Nonlinear Theory and its Applications,
  NOLTA’95, 1 (1995), pp.~91--96.

\bibitem{HOPM}
{\sc L.~De~Lathauwer, B.~De~Moor, and J.~Vandewalle}, {\em On the best rank-1
  and rank-(r 1, r 2,..., rn) approximation of higher-order tensors}, SIAM
  journal on Matrix Analysis and Applications, 21 (2000), pp.~1324--1342.

\bibitem{convergence}
{\sc A.~Falco and A.~Nouy}, {\em A proper generalized decomposition for the
  solution of elliptic problems in abstract form by using a functional
  eckart--young approach}, Journal of Mathematical Analysis and Applications,
  376 (2011), pp.~469--480.

\bibitem{survey}
{\sc L.~Grasedyck, D.~Kressner, and C.~Tobler}, {\em A literature survey of
  low-rank tensor approximation techniques}, GAMM-Mitteilungen, 36 (2013),
  pp.~53--78.

\bibitem{NP1}
{\sc J.~Haastad}, {\em Tensor rank is np-complete}, in International Colloquium
  on Automata, Languages, and Programming, Springer, 1989, pp.~451--460.

\bibitem{param3}
{\sc F.~Hamilton, T.~Berry, N.~Peixoto, and T.~Sauer}, {\em Real-time tracking
  of neuronal network structure using data assimilation}, Physical Review E, 88
  (2013), p.~052715.

\bibitem{modelfree}
{\sc F.~Hamilton, T.~Berry, and T.~Sauer}, {\em Ensemble kalman filtering
  without a model}, Physical Review X, 6 (2016), p.~011021.

\bibitem{harlim2010catastrophic}
{\sc J.~Harlim, A.~J. Majda, et~al.}, {\em Catastrophic filter divergence in
  filtering nonlinear dissipative systems}, Communications in Mathematical
  Sciences, 8 (2010), pp.~27--43.

\bibitem{NP2}
{\sc C.~J. Hillar and L.-H. Lim}, {\em Most tensor problems are np-hard},
  Journal of the ACM (JACM), 60 (2013), pp.~1--39.

\bibitem{julier1998skewed}
{\sc S.~J. Julier}, {\em Skewed approach to filtering}, in Signal and Data
  Processing of Small Targets 1998, vol.~3373, International Society for Optics
  and Photonics, 1998, pp.~271--282.

\bibitem{julier2002scaled}
{\sc S.~J. Julier}, {\em The scaled unscented transformation}, in Proceedings
  of the 2002 American Control Conference (IEEE Cat. No. CH37301), vol.~6,
  IEEE, 2002, pp.~4555--4559.

\bibitem{julier}
{\sc S.~J. Julier and J.~K. Uhlmann}, {\em A general method for approximating
  nonlinear transformations of probability distributions}, tech. report,
  Technical report, Robotics Research Group, Department of Engineering Science,
  1996.

\bibitem{julier1997}
{\sc S.~J. Julier and J.~K. Uhlmann}, {\em New extension of the kalman filter
  to nonlinear systems}, in Signal processing, sensor fusion, and target
  recognition VI, vol.~3068, International Society for Optics and Photonics,
  1997, pp.~182--193.

\bibitem{filter1}
{\sc S.~J. Julier and J.~K. Uhlmann}, {\em Unscented filtering and nonlinear
  estimation}, Proceedings of the IEEE, 92 (2004), pp.~401--422.

\bibitem{julier2004unscented}
{\sc S.~J. Julier and J.~K. Uhlmann}, {\em Unscented filtering and nonlinear
  estimation}, Proceedings of the IEEE, 92 (2004), pp.~401--422.

\bibitem{julier1995}
{\sc S.~J. Julier, J.~K. Uhlmann, and H.~F. Durrant-Whyte}, {\em A new approach
  for filtering nonlinear systems}, in Proceedings of 1995 American Control
  Conference-ACC'95, vol.~3, IEEE, 1995, pp.~1628--1632.

\bibitem{EigApprox}
{\sc E.~Kofidis and P.~A. Regalia}, {\em On the best rank-1 approximation of
  higher-order supersymmetric tensors}, SIAM Journal on Matrix Analysis and
  Applications, 23 (2002), pp.~863--884.

\bibitem{subtract2}
{\sc T.~G. Kolda}, {\em A counterexample to the possibility of an extension of
  the eckart--young low-rank approximation theorem for the orthogonal rank
  tensor decomposition}, SIAM Journal on Matrix Analysis and Applications, 24
  (2003), pp.~762--767.

\bibitem{tensor1}
{\sc T.~G. Kolda and B.~W. Bader}, {\em Tensor decompositions and
  applications}, SIAM review, 51 (2009), pp.~455--500.

\bibitem{lorenz}
{\sc E.~N. Lorenz}, {\em Deterministic nonperiodic flow}, Journal of the
  atmospheric sciences, 20 (1963), pp.~130--141.

\bibitem{HOPMrevisited}
{\sc P.~A. Regalia and E.~Kofidis}, {\em The higher-order power method
  revisited: convergence proofs and effective initialization}, in 2000 IEEE
  International Conference on Acoustics, Speech, and Signal Processing.
  Proceedings (Cat. No. 00CH37100), vol.~5, IEEE, 2000, pp.~2709--2712.

\bibitem{smoother1}
{\sc S.~S{\"a}rkk{\"a}}, {\em Unscented rauch--tung--striebel smoother}, IEEE
  transactions on automatic control, 53 (2008), pp.~845--849.

\bibitem{subtract1}
{\sc A.~Stegeman and P.~Comon}, {\em Subtracting a best rank-1 approximation
  may increase tensor rank}, Linear Algebra and its Applications, 433 (2010),
  pp.~1276--1300.

\bibitem{filter2}
{\sc R.~Van Der~Merwe, A.~Doucet, N.~De~Freitas, and E.~A. Wan}, {\em The
  unscented particle filter}, in Advances in neural information processing
  systems, 2001, pp.~584--590.

\bibitem{param1}
{\sc E.~A. Wan and R.~Van Der~Merwe}, {\em The unscented kalman filter for
  nonlinear estimation}, in Proceedings of the IEEE 2000 Adaptive Systems for
  Signal Processing, Communications, and Control Symposium (Cat. No. 00EX373),
  Ieee, 2000, pp.~153--158.

\bibitem{param2}
{\sc E.~A. Wan, R.~Van Der~Merwe, and A.~T. Nelson}, {\em Dual estimation and
  the unscented transformation}, in Advances in neural information processing
  systems, 2000, pp.~666--672.

\end{thebibliography}

\end{document}